\documentclass[print,wide,eng]{draft}

\def\cM{{\cal M}}

\title{\bf Controllability of Stochastic Evolution Equations via Samplewise Deterministic Controllability of Random-Coefficient Evolution Equations}

\author{Qi L\"u\thanks{School of Mathematics, Sichuan University, Chengdu,
610064, China. E-mail: \texttt{lu@scu.edu.cn}. This author is
supported by the National Natural Science Foundation of China under grant
12625102.}, \quad Yu Wang\thanks{School of Mathematics, Southwest Jiaotong
University, Chengdu 611756, China. E-mail: \texttt{yuwangmath@163.com}. This
author is supported by the National Natural Science Foundation of China under
grant 12401589 and by Sichuan Science and Technology Program under grant
2026NSFSC0777.}}

\date{}

\hypersetup{
pdftitle={Controllability of Stochastic Evolution Equations via Samplewise Deterministic Controllability of Random-Coefficient Evolution Equations},
pdfsubject={Transfer of deterministic endpoint relations to forward stochastic evolution systems},
pdfkeywords={endpoint reachability, stochastic evolution equation, deterministic-to-stochastic lifting, exact controllability, null controllability, approximate controllability, predictable projection, martingale representation, stochastic heat equation, stochastic wave equation, convex constraints}
}

\begin{document}
\maketitle

\begin{abstract}
We study the controllability of stochastic evolution equations by reducing the problem to the samplewise deterministic controllability of random-coefficient equations. Treating the sample point as a parameter, we use deterministic controllability theory to select measurable, square-integrable controlled trajectories with prescribed endpoints. Because such pathwise deterministic controls may depend on future information and fail to be adapted, we construct an adapted lifting that converts these trajectories into admissible controls for the stochastic equation, preserving their endpoints and estimates. Predictable projection acts on the state and drift control, while the predictable representation property and a uniformly bounded predictable right inverse of the diffusion control operator realize the martingale correction. We prove quantitative lifting estimates for deterministic
$C_0$-semigroup generators with bounded predictable perturbations
on their fixed domains. The lifting transfers deterministic exact, null, and mean-square approximate reachability to the stochastic system and preserves fixed closed convex state and drift-control constraints.

Applications include null controllability for stochastic parabolic equations, mesh-independent exact null controllability for semi-discrete stochastic parabolic equations, constrained controllability to terminal random targets for stochastic parabolic equations, and exact controllability for stochastic wave equations.  
\end{abstract}



\section{Introduction}\label{sec:introduction}

We study the controllability of stochastic evolution equations. Let $T>0$, let $(\Omega,\mathcal F,\mathbb F,\mathbb P)$ be a complete filtered probability space satisfying the usual conditions, where $\mathbb F=\{\mathcal F_t\}_{0\leq t\leq T}$, and let $W(\cdot)$ be a cylindrical Wiener process on a separable Hilbert space $\mathfrak K$. 
Throughout, all Hilbert spaces are assumed to be separable.
We assume that $\mathbb F$ has the predictable representation property
with respect to $W$: every square-integrable real-valued
$\mathbb F$-martingale can be represented as its initial value plus
a stochastic integral with respect to $W$.
This property holds, in particular, for the usual augmentation of
the filtration generated by $W$.

For a Banach space $E$ and a sub-$\sigma$-field $\mathcal G\subset\mathcal F$, let $L^2_{\mathcal G}(\Omega;E)$ and $L^\infty_{\mathcal G}(\Omega;E)$ denote the square-integrable and essentially bounded strongly $\mathcal G$-measurable random variables. For $0\leq s<r\leq T$, let $L^2_{\mathbb F}(s,r;E)$ be the space of predictable processes $X$ with
$\mathbb E\int_s^r\lVert X(t)\rVert_E^2\,d t<\infty$ and let $L^\infty_{\mathbb F}(s,r;E)$ be the space of predictable processes $X$ with
${\rm esssup}_{(t,\omega)\in (s,r)\times\Omega}\lVert X(t)\rVert_E <\infty$.
Let $L^2_{\mathbb F}(\Omega;C([s,r];E))$ be the space of adapted continuous processes with
$\mathbb E\sup_{t\in[s,r]}\lVert X(t)\rVert_E^2<\infty$.
Let $C_{\mathbb F}([s,r];L^2(\Omega;H))$ denote the space of adapted processes that are continuous from $[s,r]$ into $L^2(\Omega;H)$, equipped with the norm
$\sup_{s\leq t\leq r}(\mathbb E\lVert X(t)\rVert_H^2)^{1/2}$.
All these spaces are equipped with their usual norms.

We denote bounded linear operators  ({\rm resp}. Hilbert--Schmidt operators) by $\mathcal L(E_1;E_2)$ ({\rm resp}. $\mathcal L_2(E_1;E_2)$), abbreviate $\mathcal L(E;E)$ ({\rm resp}. $\mathcal L_2(E;E)$) to $\mathcal L(E)$ ({\rm resp}. $\mathcal L_2(E)$), and write $*$ for the Hilbert adjoint. For operator families, $\lVert\cdot\rVert_\infty$ is the essential supremum of the operator norm over $(\omega,t)$. On any interval $(s,r)$, $\Pi_E$ denotes the orthogonal projection from $L^2(\Omega\times(s,r);E)$ onto $L^2_{\mathbb F}(s,r;E)$ when $E$ is a separable Hilbert space.

Let $H,U_1,U_2$ be separable Hilbert spaces. Let $\mathcal A_0:\operatorname{Dom}(\mathcal A_0)\subset H\to H$ generate a $C_0$-semigroup $\{S_0(t)\}_{t\geq0}$ on $H$, and set $\mathcal A(t,\omega)=\mathcal A_0+N(t,\omega)$ on $\operatorname{Dom}(\mathcal A_0)$. The operator families
$N:\Omega\times(0,T)\to\mathcal L(H)$,
$B:\Omega\times(0,T)\to\mathcal L(U_1;H)$,
$C:\Omega\times(0,T)\to\mathcal L(H;\mathcal L_2(\mathfrak K;H))$,
and $D:\Omega\times(0,T)\to\mathcal L(U_2;\mathcal L_2(\mathfrak K;H))$
are strongly predictable\footnote{Strong predictability of an operator family means that its application to each fixed vector yields a predictable process.} and uniformly bounded.

We consider the control system:
\begin{align}\label{eq:intro-system}
\begin{cases}
dX=(\mathcal A X+Bu+f)\,dt+(CX+Dv)\,dW(t) & \mbox{ in } (0,T], \\[2pt]
X(0)=X_0.
\end{cases}
\end{align}
Here  
$X_0\in L^2_{\mathcal F_0}(\Omega;H)$, $f\in L^2_{\mathbb F}(0,T;H)$ and the controls satisfy
$u\in L^2_{\mathbb F}(0,T;U_1)$ and
$v\in L^2_{\mathbb F}(0,T;U_2)$. By the classical well-posedness result for stochastic evolution equations, the control system \eqref{eq:intro-system} has a unique solution $X(\cdot)\in C_{\mathbb F}([0,T];L^2(\Omega;H))$ and the following estimate holds:
\begin{align}\label{eq:abstract-mild-energy1}
\lVert X\rVert_{C_{\mathbb F}([0,T];L^2(\Omega;H))}^2
\leq\mathcal C\Bigl(
\mathbb E\lVert X_0\rVert_H^2
+\mathbb E\int_0^T
\bigl(\lVert f(t)\rVert_H^2+\lVert u(t)\rVert_{U_1}^2+\lVert v(t)\rVert_{U_2}^2\bigr)\,dt
\Bigr),
\end{align}
where $\mathcal C$ depends only on $T$,
$\sup\limits_{0\leq t\leq T}\lVert S_0(t)\rVert$, and the uniform bounds of
$N,B,C,D$  (see \cite[Theorem~3.14]{Lue2021a}).

The goal is to construct adapted controls that reach a prescribed
terminal state, or to establish exact, null, and mean-square
approximate controllability for the control system \eqref{eq:intro-system}. 
Our basic strategy is to reduce the controllability problem for the stochastic evolution equation \eqref{eq:intro-system} to the samplewise deterministic controllability of a random-coefficient evolution equation. More precisely, for a fixed sample point $\omega$, the drift part of \eqref{eq:intro-system} becomes a deterministic evolution equation whose coefficients depend on the parameter $\omega$:
\begin{equation*}
\begin{cases}\partial_t\widetilde X=\mathcal A \widetilde X+B \widetilde u+f & \mbox{ in } (0,T], \\[2pt] \widetilde X(0)=X_0,\quad \widetilde X(T)=\xi(\omega).
\end{cases}
\end{equation*}
In this formulation, the sample point is treated as a parameter, while randomness enters through the coefficients and the endpoint data. If, for almost every sample point, deterministic controllability theory can select a controlled trajectory $(\widetilde X(\omega),\widetilde u(\omega))$ with the prescribed endpoints and suitable regularity, then one obtains samplewise deterministic controls for the random-coefficient evolution equation. The remaining problem is to convert these pathwise deterministic controls into adapted controls for the original stochastic equation \eqref{eq:intro-system}. This conversion is the central construction of the paper.

The main difficulty in implementing the above reduction is that a pathwise deterministic controlled pair $(\widetilde X,\widetilde u)$ may depend on future information and hence fail to be adapted, even though it satisfies the deterministic endpoint relation for each fixed sample point. Deterministic constructions often use the terminal condition $\xi(\omega)$ or future coefficients sample by sample, so the resulting control at time $t$ may already depend on information available only at $\mathcal F_T$ rather than $\mathcal F_t$. Such a control is not admissible for the stochastic equation \eqref{eq:intro-system}. Therefore, the samplewise deterministic reduction must be supplemented by an adaptation step.

We overcome this difficulty by constructing an \emph{adapted lifting}. Given a measurable, square-integrable pathwise deterministic controlled pair $(\widetilde X,\widetilde u)$, we construct an adapted pair $(X,u)$ that preserves the same endpoints and realizes the stochastic correction through the diffusion control. Predictable projection acts on the state and drift control:
\begin{equation*}
X(t)=\mathbb E[\widetilde X(t)\mid\mathcal F_t],\qquad u=\Pi_{U_1}\widetilde u.
\end{equation*}
The projected state has a continuous adapted version, and the drift
control becomes predictable. Projecting the drift alone generally
produces a terminal residual. The predictable representation property realizes this residual as a stochastic integral and supplies the martingale integrand $\Sigma$ of the projected state.
If $D$ admits a uniformly bounded predictable right inverse $R$,
we set $v=R(\Sigma-CX)$, so that $CX+Dv=\Sigma$.
This produces admissible controls for \cref{eq:intro-system}
while preserving the prescribed endpoint.

Our main result, \cref{thm:abstract-lifting}, establishes this lifting
on the fixed domain of the generator $\mathcal A_0$ and proves
contraction estimates for the state and drift control together
with an estimate for the diffusion control. If the deterministic trajectory satisfies fixed closed convex state and control constraints, projection preserves these constraints. The construction transfers deterministic exact, null, and mean-square approximate reachability; see \cref{cor:endpoint-transfer,cor:controllability-transfer}. For the internally controlled wave equation, a change of variables reduces the lifting step to the same theorem with a zero generator; see \cref{sec:application-wave}.

The framework has two substantive requirements. First, for almost every sample point, deterministic controllability theory must provide a measurable, square-integrable controlled trajectory with the required regularity and endpoint relation. Second, the diffusion control channel must be able to realize the martingale correction produced by projection. The first requirement is supplied by deterministic PDE control theory, while the second is guaranteed by the right-inverse assumption on the diffusion control operator. The technical work of the paper is concentrated on measurable selection, trajectory regularity, lifting estimates, and applications to concrete models.

We develop four applications to the controllability of stochastic partial differential equations that illustrate how the samplewise deterministic reduction and the adapted lifting work together:
\begin{enumerate}[label=(\roman*),font=\upshape]
\item We obtain the null controllability of stochastic parabolic equations with bounded
random zeroth-order coefficients, which was first established in 
\cite{TangZhang2009} by a Carleman estimate for backward stochastic parabolic equations.
\item For one-dimensional stochastic semi-discrete parabolic equations on sufficiently fine uniform grids, we prove exact null controllability with a control cost independent of the mesh; see \cref{sd:thm:stochastic-exact}.  
\item For sufficiently long times, when the drift control is subject to pointwise box constraints and the initial data are bounded and nonnegative, we reach terminal random reference states while preserving nonnegativity; see \cref{thm:positive-main}. The reference data are uniformly bounded and terminally measurable, and their controls have a positive margin.  
\item For stochastic wave equations with a bounded random potential, we obtain exact controllability in the energy space for every $T>T_0$, where the explicit geometric threshold $T_0$ is given in \cref{eq:wave-control-time}; see \cref{thm:wave-exact}. The diffusion controls act on the whole domain in both components of the first-order system, and the lifting requires no additional control time.
\end{enumerate}

\medskip 

Controllability problems for stochastic parabolic equations have attracted increasing attention in recent years. There are generally two main approaches to these problems. The first is the spectral inequality and Lebeau--Robbiano strategy \cite{Lu2011,LiuWuYangZhong2026,YangZhong2016}. Its advantage is that a single localized drift control suffices for null controllability, but the cited results assume spatially independent multiplicative noise coefficients. The second is the global Carleman estimate, which applies to equations with spatially dependent coefficients. Using this method, Barbu, R\u{a}\c{s}canu, and Tessitore \cite{BarbuRascanuTessitore2003} established early estimates and control results for stochastic parabolic equations, and Tang and Zhang \cite{TangZhang2009} developed forward and backward null-controllability results with random coefficients. The method has since been applied to fourth-order stochastic parabolic equations \cite{GaoChenLi2015,Lue2022a}, semilinear stochastic parabolic equations \cite{HernandezSantamaria2023}, equations with state-gradient terms \cite{Zhang2024c}, and one-dimensional spatially discretized stochastic parabolic equations \cite{Zhao2025,Lecaros2025}. Its advantage is applicability to general coefficients, but the cost is the need for two controls, one in the drift term and one in the diffusion term.

In this paper, we study controllability for stochastic parabolic and wave equations using the adapted lifting method described above. For parabolic equations, this method not only recovers known results in the literature (\cref{thm:main}) but also yields new results (\cref{sd:thm:stochastic-exact,thm:positive-main}). For the wave equation, the internal observability estimate of Fu, Liao, and L\"u \cite[Theorem~1.4 and Remark~1.7]{FuLiaoLu2023} supplies the samplewise controls for bounded time-dependent potentials with an explicit geometric time threshold. The lifting preserves this time condition and reaches arbitrary square-integrable terminal random states in the energy space.

The rest of this paper is organized as follows. 
\Cref{sec:measurable-common,sec:projection} establish measurable selection and adapted lifting. The continuous, semi-discrete, and constrained parabolic equations are treated in \cref{sec:application-parabolic,sec:application-semidiscrete,sec:positive-parabolic}. \Cref{sec:application-wave} proves exact controllability for the internally controlled stochastic wave equation. The semi-discrete spectral and determinant estimates are proved in \cref{app:semidiscrete-technical}.

\section{Framework: Measurable Selection and Adapted Lifting}\label{sec:framework}

Our framework consists of two main steps. First, treating the sample point as a parameter, we select a measurable, square-integrable deterministic controlled trajectory. Second, since these pathwise trajectories may depend on future information and fail to be adapted, we construct an adapted lifting to obtain admissible controls for the stochastic evolution equation.

We first specify an auxiliary control system:
\begin{align}\label{eq:intro-system1}
\begin{cases}
dX=(\mathcal A X+Bu+f)\,dt+(CX+Dv)\,dW(t) & \mbox{ in } (s,T], \\[2pt]
X(s)=X_s,
\end{cases}
\end{align}
where $s\in [0,T)$ and  $X_s\in L^2_{\mathcal F_s}(\Omega;H)$.  Similar to \eqref{eq:intro-system}, \eqref{eq:intro-system1} has a unique solution $X(\cdot)\in C_{\mathbb F}([s,T];L^2(\Omega;H))$ and the following estimate holds:
\begin{align}\label{eq:abstract-mild-energy}
\lVert X\rVert_{C_{\mathbb F}([s,T];L^2(\Omega;H))}^2
\leq\mathcal C\Bigl(
\mathbb E\lVert X_s\rVert_H^2
+\mathbb E\int_s^T
\bigl(\lVert f(t)\rVert_H^2+\lVert u(t)\rVert_{U_1}^2+\lVert v(t)\rVert_{U_2}^2\bigr)\,dt
\Bigr),
\end{align}
where $\mathcal C$ depends only on $T$,
$\sup\limits_{0\leq t\leq T}\lVert S_0(t)\rVert$, and the uniform bounds of
$N,B,C,D$  (see \cite[Theorem~3.14]{Lue2021a}).

\subsection{Measurable selection of deterministic controlled trajectories}\label{sec:measurable-common}

To apply a pathwise controllability result, we must select its controls measurably and obtain square-integrable controlled trajectories. For a fixed sample point $\omega$, we consider the corresponding deterministic system with frozen coefficients. Its terminal state can be written as
\begin{equation*}
X_T(\omega)=X_T^0(\omega)+\Lambda_T(\omega)u,
\end{equation*}
where $X_T^0(\omega)\in H$ is the zero-input terminal state and
$\Lambda_T(\omega)\in\mathcal L(L^2(0,T;U_1);H)$ is the bounded linear operator mapping the drift control $u$ to its contribution to the terminal state. 
For almost every $\omega$, reaching the target $\xi(\omega)\in H$ amounts to finding a control $u(\omega)\in L^2(0,T;U_1)$ that solves the deterministic equation
\begin{equation*}
\Lambda_T(\omega)u(\omega)=g(\omega):=\xi(\omega)-X_T^0(\omega).
\end{equation*}
We first treat the case where $\Lambda_T$ is a random operator and the controls are unconstrained, using Tikhonov regularization to select minimum-norm controls. We then consider a fixed deterministic operator subject to convex control constraints. In both cases, deterministic solution estimates will provide the trajectory regularity required for the lifting.

To continue, we recall that an operator family $L:\Omega\to\mathcal L(U_0;E_0)$ between Hilbert spaces
is called {\it strong-operator $\mathcal G$-measurable} if
$\omega\mapsto L(\omega)u$ is strongly $\mathcal G$-measurable as an
$E_0$-valued map for every fixed $u\in U_0$.

Measurability of generalized inverses is classical
\cite{NashedSalehi1973}. The following proposition gives the selection
result used here. Tikhonov regularization produces measurable approximate
controls and converges to the minimum-norm control for attainable targets,
without a closed-range assumption. 
\begin{proposition}
\label{thm:random-Tikhonov}
Let $U_0,E_0$ be separable Hilbert spaces, let $\mathcal G$ be a sub-$\sigma$-field of $\mathcal{F}$, and let $L:\Omega\to\mathcal L(U_0;E_0)$ be strong-operator $\mathcal G$-measurable with $\lVert L(\omega)\rVert_{\mathcal L(U_0;E_0)}\leq R$ almost surely for some $R>0$. For $g\in L^2_{\mathcal G}(\Omega;E_0)$ and integers $n\geq1$, define\vspace{-2mm}
\begin{equation}\label{eq:common-Tikhonov}
u_n:=L^*(LL^*+n^{-1}\operatorname{Id}_{E_0})^{-1}g.\vspace{-2mm}
\end{equation}
Then $u_n\in L^2_{\mathcal G}(\Omega;U_0)$ and\vspace{-2mm}
\begin{equation}\label{eq:dense-range-control-bound}
\bigl\lVert u_n(\omega)\bigr\rVert_{U_0}\leq\frac{\sqrt n}{2}\bigl\lVert g(\omega)\bigr\rVert_{E_0}.\vspace{-2mm}
\end{equation}
\begin{enumerate}[label=(\roman*),font=\upshape]
\item If $g(\omega)\in\operatorname{Ran} L(\omega)$ almost surely, then,
for almost every $\omega$,
$u_n(\omega)\to u_\dagger(\omega)$ in $U_0$, where
$u_\dagger(\omega)$ is the minimum-norm solution:\vspace{-2mm}
\begin{equation*}
u_\dagger(\omega)=\operatorname*{argmin}\bigl\{\lVert u\rVert_{U_0}:u\in U_0,\ L(\omega)u=g(\omega)\bigr\}.\vspace{-2mm}
\end{equation*}
The map $u_\dagger$ has a strongly $\mathcal G$-measurable version. If
$\bigl\lVert u_\dagger(\omega)\bigr\rVert_{U_0}\leq q(\omega)$ for some
$q\in L^2_{\mathcal G}(\Omega)$, then $u_\dagger$ is square integrable.
\item In general,
$L(\omega)u_n(\omega)\to P_{\overline{\operatorname{Ran}L(\omega)}}g(\omega)$
in $E_0$ for almost every $\omega$, where
$P_{\overline{\operatorname{Ran}L(\omega)}}$ is the orthogonal projection
onto the closure of $\operatorname{Ran}L(\omega)$.
The convergence also holds in $L^2_{\mathcal G}(\Omega;E_0)$.
In particular, if $L(\omega)$ has dense range almost surely, then
$Lu_n\to g$ in $L^2_{\mathcal G}(\Omega;E_0)$.
\end{enumerate}
\end{proposition}

\begin{proof} The proof is divided into three steps. 

\textbf{Step 1. Construct measurable regularized controls.}
For each fixed $\omega$,  define a functional $J_\omega: U_0\to \mathbb R$ as follows:
\begin{equation*}
J_\omega(u):=\lVert L(\omega)u-g(\omega)\rVert_{E_0}^2+n^{-1}\lVert u\rVert_{U_0}^2,\qquad u\in U_0.
\end{equation*}
Strict convexity and coercivity of $J_\omega(\cdot)$ guarantee that $J_\omega(\cdot)$
has a unique minimizer over $u\in U_0$. The first-order optimality condition for this minimizer yields
\begin{equation*}
L(\omega)^*(L(\omega)u-g(\omega))+n^{-1}u=0.
\end{equation*}
Solving for $u$ gives \cref{eq:common-Tikhonov}, which therefore defines the unique minimizer of $J_\omega$.

To prove measurability, note that for each $y\in E_0$, the identity $\langle L(\omega)^*y,u\rangle_{U_0}=\langle y,L(\omega)u\rangle_{E_0}$ for $u\in U_0$, together with the Pettis measurability theorem (\cite[Chapter 2, Section 1, Theorem 2]{DiestelUhl1977}), shows that $L^*$ is strong-operator measurable. 

To show that applying $L$ or $L^*$ to a strongly measurable random vector preserves measurability, let $\eta:\Omega\to U_0$ be strongly measurable. By the Pettis  measurability theorem (\cite[Chapter 2, Section 1, Theorem 2]{DiestelUhl1977}), $\eta$ can be approximated almost everywhere by a sequence of simple random vectors $\eta_m=\sum\limits_{k=1}^{N_m}\mathbf{1}_{A_{m,k}}u_{m,k}$, where $A_{m,k}\in\mathcal{G}$ and $u_{m,k}\in U_0$. Since $L$ is strong-operator measurable, for each fixed $u_{m,k}$, the map $\omega\mapsto L(\omega)u_{m,k}$ is strongly measurable. Consequently, the map $\omega\mapsto L(\omega)\eta_m(\omega)=\sum\limits_{k=1}^{N_m}\mathbf{1}_{A_{m,k}}L(\omega)u_{m,k}$ is strongly measurable. Because $L(\omega)$ is uniformly bounded in operator norm and $\eta_m\to \eta$ almost everywhere, we have $$\lVert L(\omega)\eta_m(\omega)-L(\omega)\eta(\omega)\rVert_{E_0}\leq R\lVert \eta_m(\omega)-\eta(\omega)\rVert_{U_0}\to0 \mbox{  almost everywhere.}$$ The pointwise limit of strongly measurable functions is strongly measurable, so $\omega\mapsto L(\omega)\eta(\omega)$ is strongly measurable. An analogous argument holds for $L^*$ using its uniform boundedness.
Applying this to the strongly measurable random vector $\eta(\omega)=L(\omega)^*y$ for each fixed $y\in E_0$ shows that $LL^*$ is strong-operator measurable. Iterating this argument shows that all powers of $LL^*$ are also strong-operator measurable.

Write $I=\operatorname{Id}_{E_0}$. By the assumed bound on $L$, the series
\begin{equation*}
(LL^*+\varepsilon I)^{-1}
=\frac1{R^2+\varepsilon}
\sum_{j=0}^\infty
\biggl(\frac{R^2I-LL^*}{R^2+\varepsilon}\biggr)^j
\end{equation*}
converges in operator norm uniformly on a set of full probability for each fixed $\varepsilon>0$. 

Since $LL^*$ is strong-operator measurable, each term in the series is strong-operator measurable. The uniform convergence on a full-probability set implies that the limit $(LL^*+\varepsilon I)^{-1}$ is also strong-operator measurable (as it is the uniform limit of measurable operator-valued functions). Applying this operator to the strongly measurable vector $g(\omega)$ preserves measurability, so the map $\omega\mapsto (L(\omega)L(\omega)^*+\varepsilon I)^{-1}g(\omega)$ is strongly measurable. Finally, composing with the strong-operator measurable $L^*$ yields the strong measurability of $u_n$.   
To estimate its norm, fix $\omega$ and suppress it in the notation. Since $LL^*$ is nonnegative and self-adjoint, we can apply the functional calculus for bounded self-adjoint operators (the spectral theorem). Let $T = LL^*$ and define the scalar function $f(\lambda) = \frac{\lambda}{(\lambda + n^{-1})^2}$ for $\lambda \ge 0$. By \cref{eq:common-Tikhonov}, we can express the squared norm of $u_n$ as
\begin{align*}
\lVert u_n\rVert_{U_0}^2
&= \langle L^*(T + n^{-1}I)^{-1}g, L^*(T + n^{-1}I)^{-1}g \rangle_{U_0} \\
&= \langle (T + n^{-1}I)^{-1} T (T + n^{-1}I)^{-1} g, g \rangle_{E_0} \\
&= \langle f(T)g, g \rangle_{E_0}.
\end{align*}
By the spectral theorem, there exists a projection-valued measure $E$ such that $f(T) = \int_{[0, \lVert T \rVert]} f(\lambda) \, dE(\lambda)$. Consequently,
\begin{align*}
\lVert u_n\rVert_{U_0}^2
&= \int_{[0, \lVert T \rVert]} \frac{\lambda}{(\lambda + n^{-1})^2} \, d\langle E(\lambda)g, g \rangle_{E_0} \\
&\leq \biggl[\sup_{\lambda \geq 0} \frac{\lambda}{(\lambda + n^{-1})^2}\biggr] \int_{[0, \lVert T \rVert]} \, d\langle E(\lambda)g, g \rangle_{E_0} \\
&= \biggl[\sup_{\lambda \geq 0} \frac{\lambda}{(\lambda + n^{-1})^2}\biggr] \lVert g \rVert_{E_0}^2.
\end{align*}
A direct differentiation of the scalar function $f(\lambda)$ shows that its maximum over $\lambda \geq 0$ is attained at $\lambda = n^{-1}$, and the maximum value is $\frac{n}{4}$. Therefore,
\begin{equation*}
\lVert u_n\rVert_{U_0}^2 \leq \frac{n}{4}\lVert g \rVert_{E_0}^2,
\end{equation*}
which proves \cref{eq:dense-range-control-bound} and, since $g\in L^2_{\mathcal G}(\Omega;E_0)$, also yields $u_n\in L^2_{\mathcal G}(\Omega;U_0)$.

\smallskip

\textbf{Step 2. Pass to the minimum-norm control.}

We first clarify the geometry of the solution set. If $g(\omega)\in\operatorname{Ran}L(\omega)$, the set of solutions $\mathcal S_\omega:=\{u\in U_0:L(\omega)u=g(\omega)\}$ is a nonempty closed affine subspace of $U_0$. By the Hilbert projection theorem (\cite[Theorem 5.2]{Brezis2011}), $\mathcal S_\omega$ contains a unique element of minimal norm, denoted $u_\dagger(\omega)$. Since $\mathcal S_\omega$ is an affine translate of the closed subspace $\operatorname{Ker}L(\omega)$, this minimum-norm element is characterized by the orthogonality condition $u_\dagger(\omega)\perp\operatorname{Ker}L(\omega)$.

Substituting $g=Lu_\dagger$ into \cref{eq:common-Tikhonov} and using the identity $L^*(LL^*+n^{-1}I)^{-1} = (L^*L+n^{-1}\operatorname{Id}_{U_0})^{-1}L^*$, we obtain
\begin{equation*}
u_n  = L^*(LL^*+n^{-1}I)^{-1}Lu_\dagger 
= L^*L(L^*L+n^{-1}\operatorname{Id}_{U_0})^{-1}u_\dagger.
\end{equation*}

Next, we examine the spectral properties of the self-adjoint nonnegative operator $L^*L$. By the spectral decomposition theorem  (\cite[Theorem~12.22]{Rudin1991}), there exists a projection-valued measure $E$ such that
\begin{equation*}
u_n = \int_{[0,\lVert L\rVert^2]} \frac{\lambda}{\lambda+n^{-1}}\,d E_\lambda u_\dagger.
\end{equation*}
The spectral multipliers $f_n(\lambda):=\lambda/(\lambda+n^{-1})$ satisfy $0\leq f_n(\lambda)\leq 1$ for all $\lambda\geq0$, and $\lim\limits_{n\to\infty}f_n(\lambda)=1$ for every $\lambda>0$.

Since $u_\dagger\perp\operatorname{Ker}L$ and $\operatorname{Ker}(L^*L)=\operatorname{Ker}L$, the spectral measure of $u_\dagger$ assigns zero mass to $\{0\}$. Thus, $f_n(\lambda)\to 1$ holds $E$-almost everywhere with respect to $u_\dagger$. By the dominated convergence theorem for spectral integrals, the boundedness of $f_n$ and the pointwise convergence guarantee that
\begin{equation*}
u_n \longrightarrow \int_{[0,\lVert L\rVert^2]} 1\,d E_\lambda u_\dagger = u_\dagger \quad \text{in } U_0 \text{ as } n\to\infty.
\end{equation*}
Finally, since $u_\dagger$ is strongly measurable and the assumed majorant controls its norm, $u_\dagger$ is square integrable. The limit of the strongly measurable sequence $u_n$ is also strongly measurable.

\smallskip

\textbf{Step 3. Identify the approximate endpoint.}

For arbitrary $g$, we analyze the limit of $Lu_n$. First, we introduce a scalar function $\phi_n(\lambda) := \frac{\lambda}{\lambda + n^{-1}}$ for $\lambda \geq 0$. Using the functional calculus for the nonnegative self-adjoint operator $LL^*$, we can express
\begin{equation*}
LL^*(LL^* + n^{-1}I)^{-1} = \phi_n(LL^*).
\end{equation*}
By the spectral decomposition theorem  (\cite[Theorem~12.22]{Rudin1991}), there exists a projection-valued measure $E$ such that
\begin{equation*}
Lu_n = \phi_n(LL^*)g = \int_{\sigma(LL^*)} \frac{\lambda}{\lambda+n^{-1}} \, dE_\lambda g.
\end{equation*}
The scalar functions $\phi_n$ satisfy $0 \leq \phi_n(\lambda) \leq 1$ for all $\lambda \geq 0$, and $\lim\limits_{n\to\infty} \phi_n(\lambda) = \chi_{(0,\infty)}(\lambda)$ for every $\lambda \geq 0$. Since $\overline{\operatorname{Ran} L}=(\operatorname{Ker}L^*)^\perp=\overline{\operatorname{Ran}LL^*}$, the projection $\chi_{(0,\infty)}(LL^*)$ corresponds exactly to the orthogonal projection onto the closure of $\operatorname{Ran} L$, denoted by $P_{\overline{\operatorname{Ran} L}}$. By the dominated convergence theorem for spectral integrals, the pointwise convergence of $\phi_n$ and the uniform bound $\phi_n \leq 1$ yield the strong convergence
\begin{equation*}
Lu_n = \phi_n(LL^*)g \longrightarrow \chi_{(0,\infty)}(LL^*)g = P_{\overline{\operatorname{Ran} L}}g \qquad\mbox{ in $E_0$ for almost every $\omega$.}
\end{equation*} 
To obtain the mean-square convergence, we note that $\lVert Lu_n \rVert_{E_0} = \lVert \phi_n(LL^*)g \rVert_{E_0} \leq \lVert g \rVert_{E_0}$. Since $g \in L^2_{\mathcal G}(\Omega;E_0)$, the uniform bound provides a square-integrable majorant. Therefore, by the Lebesgue dominated convergence theorem, we conclude that
\begin{equation*}
Lu_n \longrightarrow P_{\overline{\operatorname{Ran} L}}g \quad \text{in } L^2_{\mathcal G}(\Omega;E_0) \text{ as } n\to\infty.
\end{equation*}
\end{proof}
\begin{remark}
\label{rem:selector-integrability-gap}
A measurable selection need not be square integrable, even when the target is. For example, on $\Omega=(0,1)$ equipped with the Lebesgue measure, take $U_0=E_0=\mathbb{R}$, $L(\omega)u=\omega u$, and $g(\omega)=\sqrt\omega$. Then $L$ is uniformly bounded and measurable, $g\in L^2(\Omega)$, and $g(\omega)\in\operatorname{Ran}L(\omega)$. However, the unique minimum-norm solution is $u_\dagger(\omega)=\omega^{-1/2}$, which is not in $L^2(\Omega)$. Therefore, an additional square-integrable control bound is necessary.
\end{remark}

To enforce bounded convex control constraints, we instead minimize the control norm over a fixed nonempty, bounded, closed, and convex admissible set $K \subset U_0$. In the next result, the terminal control operator and the set $K$ are deterministic, while the target may be random. The boundedness of $K$ guarantees square integrability once the selection is measurable.
\begin{proposition}
\label{prop:constrained-selector}
Let $U_0$ and $E_0$ be separable Hilbert spaces, let
$L\in\mathcal L(U_0;E_0)$, and let $K\subset U_0$ be nonempty, bounded,
closed, and convex. Then $L(K)$ is a closed convex subset of $E_0$.
For $g\in L(K)$, let  
\begin{equation*}
\mathfrak S_K(g):=\operatorname*{argmin}\bigl\{\lVert u\rVert_{U_0}:u\in K,\ Lu=g\bigr\}.
\end{equation*}
Then $\mathfrak S_K:L(K)\to U_0$ is Borel measurable.
Consequently, if $\mathcal G$ is a sub-$\sigma$-field
of $\mathcal{F}$, $g:\Omega\to E_0$ is strongly
$\mathcal G$-measurable and $g(\omega)\in L(K)$ almost surely, then
$\mathfrak S_K(g)$ has a strongly $\mathcal G$-measurable version with
values in $K$, belongs to $L^2_{\mathcal G}(\Omega;U_0)$, and satisfies
$L\mathfrak S_K(g)=g$ almost surely.
\end{proposition}
\begin{proof}
\textbf{Step 1. Regularize the constrained problem.}
Since $U_0$ is a Hilbert space (and hence reflexive), the bounded, closed, and convex set $K$ is weakly compact by the Banach-Alaoglu theorem. Because $L$ is a bounded linear operator, it is weakly continuous; thus, its image $L(K)$ is also weakly compact. 

Being weakly compact, $L(K)$ is weakly closed and norm bounded. Moreover, the linearity of $L$ preserves the convexity of $K$. By Mazur's theorem, the norm closure of a convex set coincides with its weak closure. Hence, $L(K)$ is closed in the norm topology. Therefore, $L(K)$ is a nonempty, bounded, closed, and convex subset of $E_0$. 

Consequently, for every $g\in L(K)$, the admissible solution set $\{u\in K:Lu=g\}$ is a nonempty closed convex subset of $U_0$. By the Hilbert projection theorem, it admits a unique element of minimal norm.

For $m\geq1$, consider the regularized Tikhonov functional
\begin{equation*}
J_m(u;g) := \lVert Lu-g\rVert_{E_0}^2 + m^{-1}\lVert u\rVert_{U_0}^2, \qquad u\in K.
\end{equation*}
Since the norm is weakly lower semicontinuous and $L$ is a bounded linear operator, the functional $J_m(\cdot;g)$ is also weakly lower semicontinuous. Therefore, by the extreme value theorem in the weak topology, $J_m(\cdot;g)$ attains its minimum at some point in the weakly compact set $K$. 

The term $m^{-1}\lVert u\rVert_{U_0}^2$ is strictly convex in $u$, which makes the entire functional $J_m(\cdot;g)$ strictly convex. Since a strictly convex functional can have at most one minimizer on a convex set, the minimizer is unique. We denote it by
\begin{equation*}
u_m(g)=\operatorname*{argmin}_{u\in K}
\bigl\{\lVert Lu-g\rVert_{E_0}^2+m^{-1}\lVert u\rVert_{U_0}^2\bigr\},
\qquad g\in E_0.
\end{equation*}
To prove that $u_m(\cdot)$ is continuous, set $u_i=u_m(g_i)$ for $i=1,2$. The first-order variational inequality for $u_i$ yields
\begin{align}\label{11.10-eq1}
\operatorname{Re}\langle L^*(Lu_i-g_i)+m^{-1}u_i,w-u_i\rangle_{U_0}
\geq0,\qquad w\in K,\quad i=1,2.
\end{align}
Taking $w=u_2$ and $w=u_1$ in \eqref{11.10-eq1}, respectively, and adding the resulting inequalities gives
\begin{equation*}
\lVert L(u_1-u_2)\rVert_{E_0}^2+m^{-1}\lVert u_1-u_2\rVert_{U_0}^2
 \leq\operatorname{Re}\langle g_1-g_2,L(u_1-u_2)\rangle_{E_0}.
\end{equation*}
By the Cauchy--Schwarz inequality, this implies
$$\lVert u_m(g_1)-u_m(g_2)\rVert_{U_0}
\leq m\lVert L\rVert\lVert g_1-g_2\rVert_{E_0},$$ which yields that $u_m$ is continuous.

\smallskip

\textbf{Step 2. Pass to the constrained minimum-norm control.}

Let $g\in L(K)$ and let $u_*=\mathfrak S_K(g)$ be the unique minimum-norm solution. The minimality of $u_m(g)$ for the regularized functional yields
\begin{align}\label{eq:constrained-minimality}
\bigl\lVert Lu_m(g)-g\bigr\rVert_{E_0}^2+m^{-1}\bigl\lVert u_m(g)\bigr\rVert_{U_0}^2
\leq J_m(u_*;g) = m^{-1}\lVert u_*\rVert_{U_0}^2.
\end{align}
From \cref{eq:constrained-minimality}, we directly deduce two bounds:
\begin{enumerate}[label=(\roman*),font=\upshape]
\item $\bigl\lVert Lu_m(g)-g\bigr\rVert_{E_0}^2 \leq m^{-1}\lVert u_*\rVert_{U_0}^2$, which implies $Lu_m(g)\to g$ in $E_0$ as $m\to\infty$.
\item $m^{-1}\bigl\lVert u_m(g)\bigr\rVert_{U_0}^2 \leq m^{-1}\lVert u_*\rVert_{U_0}^2$, which implies $\bigl\lVert u_m(g)\bigr\rVert_{U_0} \leq \lVert u_*\rVert_{U_0}$. Thus, the sequence $\{u_m(g)\}_{m=1}^\infty$ is bounded uniformly in $m$.
\end{enumerate}
Since $K$ is   weakly compact, the bounded sequence $\{u_m(g)\}_{m=1}^\infty$ admits a weakly convergent subsequence $\{u_{m_k}(g)\}_{k=1}^\infty$ such that $u_{m_k}(g) \rightharpoonup \widehat{u}$ for some $\widehat{u}\in K$ as $k\to\infty$. Because $L$ is a bounded linear operator, it is weakly continuous. Hence $Lu_{m_k}(g) \rightharpoonup L\widehat{u}$ as $k\to\infty$. But we already know $Lu_m(g)\to g$ strongly (and thus weakly) as $m\to\infty$. Thus, by the uniqueness of weak limits, we must have $L\widehat{u}=g$. This shows that $\widehat{u}$ is feasible (i.e., $\widehat{u}\in \{u\in K: Lu=g\}$).

Next, by the weak lower semicontinuity of the norm, we have
\begin{equation*}
\lVert \widehat{u}\rVert_{U_0} \leq \liminf_{k\to\infty} \bigl\lVert u_{m_k}(g)\bigr\rVert_{U_0} \leq \lVert u_*\rVert_{U_0}.
\end{equation*}
Since $u_*$ is the  unique  element of minimal norm in the feasible set, the inequality $\lVert \widehat{u}\rVert_{U_0} \leq \lVert u_*\rVert_{U_0}$ forces $\widehat{u}=u_*$. This proves that $u_*$ is the only possible weak cluster point of $\{u_m(g)\}_{m=1}^\infty$. Since every subsequence has a further subsequence converging weakly to $u_*$, the entire sequence converges weakly: $u_m(g) \rightharpoonup u_*$ as $m\to\infty$.

Finally, we upgrade this weak convergence to strong convergence. By weak lower semicontinuity and the norm bound from \cref{eq:constrained-minimality},
\begin{equation*}
\lVert u_*\rVert_{U_0}
\leq \liminf_{m\to\infty}\lVert u_m(g)\rVert_{U_0}
\leq \limsup_{m\to\infty}\lVert u_m(g)\rVert_{U_0}
\leq \lVert u_*\rVert_{U_0}.
\end{equation*}
Consequently,  $\lVert u_m(g)\rVert_{U_0} \to \lVert u_*\rVert_{U_0}$  as $m\to\infty$. Since $u_m(g) \rightharpoonup u_*$ weakly  as $m\to\infty$ and their norms converge, a standard property of Hilbert spaces implies the strong convergence $u_m(g)\to u_*$ in $U_0$ as $m\to\infty$.

The pointwise limit $\mathfrak S_K$ is thus Borel measurable on the closed set $L(K)$. To see this, recall that the sequence $\{u_m(\cdot)\}_{m=1}^\infty$ consists of continuous (hence Borel measurable) maps from $E_0$ to $U_0$. Since $u_m(g) \to \mathfrak S_K(g)$ as $m\to\infty$ for every $g \in L(K)$, the pointwise limit of Borel measurable functions is Borel measurable. It can be extended by a fixed element of $K$ outside this set.  
To apply this to the random target, let $g:\Omega \to E_0$ be a strongly $\mathcal G$-measurable random variable with $g(\omega) \in L(K)$ almost surely. Since $\mathfrak S_K$ is Borel measurable on $L(K)$ and $g$ is strongly $\mathcal G$-measurable, the standard composition property of measurable functions ensures that the map $\omega \mapsto \mathfrak S_K(g(\omega))$ is strongly $\mathcal G$-measurable.

Finally, because the set $K$ is bounded, there exists a constant $C_K > 0$ such that $\lVert u \rVert_{U_0} \leq C_K$ for all $u \in K$. Therefore, the selection satisfies $\lVert \mathfrak S_K(g(\omega)) \rVert_{U_0} \leq C_K$ almost surely. Since $g$ is strongly $\mathcal G$-measurable, the composition $\mathfrak S_K(g)$ is strongly $\mathcal G$-measurable. Combining this with the uniform bound $C_K$, we conclude that $\mathfrak S_K(g) \in L^2_{\mathcal G}(\Omega; U_0)$.
\end{proof}
In the applications, we combine these selectors with deterministic solution
operators and regularity estimates to obtain measurable, square-integrable
controlled trajectories. These trajectories may depend on future information
and need not be adapted. The next subsection addresses this remaining issue.

\subsection{Adapted lifting of deterministic controlled trajectories}\label{sec:projection}

We now lift a measurable, square-integrable controlled pair
$(\widetilde X,\widetilde u)$ to an adapted pair $(X,u)$ with the same
endpoints. Predictable projection determines the drift, and the
predictable representation property supplies a stochastic integral to
compensate for the terminal residual. On the common operator
domain, we project the instantaneous drift. A right inverse
of the diffusion control operator then realizes this correction through
$CX+Dv=\Sigma$.

\subsubsection{The predictable representation property and predictable projection}

In the rest of this paper, generic constants  $\mathcal{C}$ in estimates are positive and may change from line
to line. Their relevant dependence is specified in each application.

For $0\leq s<r\leq T$ and
$\xi\in L^2_{\mathcal{F}_r}(\Omega;H)$ satisfying $\mathbb{E}[\xi\mid\mathcal{F}_s]=0$,
the predictable representation property gives a unique
$Z\in L^2_{\mathbb F}(s,r;\mathcal L_2(\mathfrak K;H))$ such that\vspace{-2mm}
\begin{align}\label{eq:shifted-representation}
\xi=\int_s^r Z(t)\,d W(t).
\end{align}
It\^o's isometry gives\vspace{-2mm}
\begin{equation*}
\mathbb{E}\int_s^r\bigl\lVert Z(t)\bigr\rVert_{\mathcal{L}_2(\mathfrak K;H)}^2\,d t
=\mathbb{E}\lVert \xi\rVert_H^2.
\end{equation*}

The following lemma identifies the semimartingale decomposition of the
conditional expectation of an $H^1$ trajectory and bounds its martingale
term. We combine classical projection theory
\cite[Chapter~VI, Section~2]{DellacherieMeyer1982} with the predictable
representation property, keeping the endpoint identities explicit.
\begin{lemma}\label{thm:causal-projection}
Let $0\leq s<r\leq T$ and
$\widetilde X\in L^2_{\mathcal F_r}(\Omega;H^1(s,r;H))$,
with $\widetilde X(s)\in L^2_{\mathcal F_s}(\Omega;H)$.
The conditional expectation process\vspace{-2mm}
\begin{align}\label{eq:conditional-path}
X(t):=\mathbb E[\widetilde X(t)\mid\mathcal F_t],
\qquad t\in[s,r],
\end{align}
admits a continuous adapted version, still denoted by $X$, belonging to
$L^2_{\mathbb F}(\Omega;C([s,r];H))$.

Moreover, there exists a unique 
$Z\in L^2_{\mathbb F}(s,r;\mathcal L_2(\mathfrak K;H))$
such that $X$ admits the semimartingale decomposition: almost surely for all $t\in[s,r]$,\vspace{-2mm}
\begin{align}\label{eq:causal-decomposition}
X(t)=\widetilde X(s)
+\int_s^t(\Pi_H\dot{\widetilde X})(\rho)\,d\rho
+\int_s^t Z(\rho)\,dW(\rho).
\end{align}
In particular,   $X(s)=\widetilde X(s)$ and
$X(r)=\widetilde X(r)$ almost surely, and $X=\Pi_H\widetilde X$ in $L^2(\Omega\times(s,r);H)$.\footnote{Recall that on any interval $(s,r)$, $\Pi_H$ denotes the orthogonal projection from $L^2(\Omega\times(s,r);H)$ onto $L^2_{\mathbb F}(s,r;H)$.}
Furthermore, the martingale integrand satisfies the estimate\vspace{-2mm}
\begin{align}\label{eq:martingale-defect-estimate}
\mathbb E\int_s^r\lVert Z(t)\rVert_{\mathcal L_2(\mathfrak K;H)}^2\,dt
\leq(r-s)\mathbb E\int_s^r
\lVert \dot{\widetilde X}(t)\rVert_H^2\,dt.
\end{align}
\end{lemma}
\begin{proof}
\textbf{Step 1. Construct the residual and apply the predictable representation property.}

Set $g=\dot{\widetilde X}$. Since $\widetilde X\in L^2_{\mathcal F_r}(\Omega;H^1(s,r;H))$, the time derivative satisfies $g\in L^2(\Omega\times(s,r);H)$. Let $\Pi_H g \in L^2_{\mathbb F}(s,r;H)$ be its orthogonal projection onto the predictable processes. Consequently, the difference $g-\Pi_H g$ lies in the orthogonal complement of the predictable subspace in $L^2(\Omega\times(s,r);H)$.

For every fixed $t\in[s,r]$ and any $\eta\in L^2_{\mathcal F_t}(\Omega;H)$, the process $\mathbf{1}_{[t,r]}(\rho)\eta$ is predictable. Since $\Pi_H$ is the orthogonal projection onto the subspace of predictable processes, the residual $g - \Pi_H g$ is orthogonal to any predictable process. Taking the test process $\mathbf{1}_{[t,r]}(\rho)\eta$ with $\eta \in L^2_{\mathcal F_t}(\Omega;H)$, we obtain \vspace{-2mm}
\begin{equation*}
\mathbb E \int_t^r \bigl\langle g(\rho)-\Pi_H g(\rho), \eta \bigr\rangle_H \,d\rho = 0.
\end{equation*}
Since this holds for every $\mathcal F_t$-measurable $\eta$, it implies the conditional expectation identity\vspace{-2mm}
\begin{align}\label{eq:projection-future-residual}
\mathbb E\Bigl[\int_t^r \bigl(g(\rho)-\Pi_H g(\rho)\bigr)\,d\rho \,\Bigm|\,\mathcal F_t\Bigr]=0.
\end{align}

Now define the terminal residual\vspace{-4mm}
\begin{equation*}
\zeta:=\int_s^r \bigl(g(\rho)-\Pi_H g(\rho)\bigr)\,d\rho.
\end{equation*}
Because $g$ and $\Pi_H g$ belong to $L^2(\Omega\times(s,r);H)$,   $\zeta$ is a well-defined square-integrable $\mathcal F_r$-measurable random variable, i.e., $\zeta \in L^2_{\mathcal F_r}(\Omega;H)$. Taking $t=s$ in \cref{eq:projection-future-residual} yields $\mathbb E[\zeta\mid\mathcal F_s]=0$. 

By the predictable representation property of the Wiener process $W$ with respect to $\mathbb F$, any square-integrable $\mathcal F_r$-measurable random variable with zero $\mathcal F_s$-conditional expectation admits a unique stochastic integral representation. Therefore, there exists a unique integrand $Z\in L^2_{\mathbb F}(s,r;\mathcal L_2(\mathfrak K;H))$ such that\vspace{-2mm}
\begin{align}\label{11.6-eq1}
\zeta=\int_s^r Z(\rho)\,d W(\rho).
\end{align}

\textbf{Step 2. Verify the conditional expectation identity.}
Define the process
\begin{align}\label{eq:def-X}
\widehat X(t):=\widetilde X(s)+\int_s^t(\Pi_H g)(\rho)\,d\rho+\int_s^t Z(\rho)\,d W(\rho), \qquad t\in[s,r].
\end{align}
By \cref{11.6-eq1,eq:projection-future-residual}, for every $t\in[s,r]$ we have
\begin{align}\label{eq:def-X1}
\int_s^tZ(\rho)\,d W(\rho)
&=\mathbb E[\zeta\mid\mathcal F_t]  =\mathbb E\Bigl[\int_s^t(g-\Pi_H g)\,d\rho\,\Bigm|\,\mathcal F_t\Bigr] + \mathbb E\Bigl[\int_t^r(g-\Pi_H g)\,d\rho\,\Bigm|\,\mathcal F_t\Bigr] \\
&=\mathbb E\Bigl[\int_s^t(g-\Pi_H g)\,d\rho\,\Bigm|\,\mathcal F_t\Bigr].
\end{align}
Since $\widetilde X(s)$ and $\int_s^t(\Pi_H g)\,d\rho$ are $\mathcal F_t$-measurable, substituting \cref{eq:def-X1} into \cref{eq:def-X} yields
\begin{equation*}
\widehat X(t) =\mathbb E\Bigl[\widetilde X(s)+\int_s^tg(\rho)\,d\rho\,\Bigm|\,\mathcal F_t\Bigr]
=\mathbb E[\widetilde X(t)\mid\mathcal F_t].
\end{equation*}
This proves that $\widehat X$ is exactly the conditional expectation process $X$ in \cref{eq:conditional-path}.

\smallskip

\textbf{Step 3. Endpoints, regularity, and uniqueness.}
At $t=s$, the integrals in \cref{eq:def-X} vanish, giving $X(s)=\widetilde X(s)$. At $t=r$, using the definition of $\zeta$, we have
\begin{align*}
X(r)&=\widetilde X(s)+\int_s^r(\Pi_H g)\,d\rho+\int_s^r Z\,d W \\
&=\widetilde X(s)+\int_s^r(\Pi_H g)\,d\rho+\int_s^r(g-\Pi_H g)\,d\rho
\\
&=\widetilde X(s)+\int_s^rg(\rho)\,d\rho=\widetilde X(r).
\end{align*}
Thus, the endpoint identities hold. The continuous version and the $L^2$ bound follow directly from \cref{eq:def-X} and Doob's maximal inequality. Uniqueness of $Z$ follows from Itô's isometry. The identity $X=\Pi_H\widetilde X$ in $L^2(\Omega\times(s,r);H)$ follows from \cref{eq:conditional-path} and the definition of the predictable projection.

\smallskip

\textbf{Step 4. Estimate the martingale defect.}
By It\^o's isometry and the Cauchy--Schwarz inequality, we obtain
\begin{align*}
\mathbb E\int_s^r\lVert Z\rVert_{\mathcal{L}_2(\mathfrak K;H)}^2\,d t
&=\mathbb E\lVert \zeta\rVert_H^2
=\mathbb E\bigl\lVert \int_s^r (g-\Pi_Hg)\,d\rho \bigr\rVert_H^2  \\
& \leq(r-s)\mathbb E\int_s^r \lVert g-\Pi_Hg\rVert_H^2\,d\rho  \\
& \leq(r-s)\mathbb E\int_s^r \lVert g\rVert_H^2\,d\rho,
\end{align*}
where the last inequality uses the contraction property of the complementary projection $I-\Pi_H$. This completes the proof.
\end{proof}

\subsubsection{Lifting on a common operator domain}\label{sec:abstract}

Fix $0\leq s<r\leq T$ and put
$\mathcal D=\operatorname{Dom}(\mathcal A_0)$,
equipped with the graph norm
$
\lVert x\rVert_{\mathcal D}^2
=\lVert x\rVert_H^2+\lVert \mathcal A_0x\rVert_H^2.
$
A fixed equivalent Hilbert norm may be used in the applications.
Recall that $\mathcal A=\mathcal A_0+N$, where $N$ is strongly
predictable and uniformly bounded in $\mathcal L(H)$.
Consequently,
$\mathcal A:\Omega\times(s,r)\to\mathcal L(\mathcal D;H)$
is strongly predictable and uniformly bounded.
The assumptions on $B$ and $C$ are as stated above.
Let $f\in L^2_{\mathbb F}(s,r;H)$. Define the space of strong trajectories
\begin{align}\label{eq:abstract-intersection-space}
\mathscr{X}(s,r)=H^1(s,r;H)\cap L^2(s,r;\mathcal D),
\end{align}
with its intersection norm. We identify an $H$-valued continuous process
with its $\mathcal D$-valued predictable representative whenever both exist.
We first take $U_2=\mathcal L_2(\mathfrak K;H)$ and
$D=\operatorname{Id}_{\mathcal L_2(\mathfrak K;H)}$.
\cref{cor:endpoint-transfer} treats diffusion control operators with a bounded
predictable right inverse. The next theorem lifts a deterministic
controlled trajectory while preserving its endpoints and fixed closed
convex constraints.
\begin{theorem}\label{thm:abstract-lifting}
Let $\widetilde X\in L^2_{\mathcal{F}_r}(\Omega;\mathscr{X}(s,r))$ with
$\widetilde X(s)\in L^2_{\mathcal{F}_s}(\Omega;H)$.
Suppose that $\widetilde u\in L^2_{\mathcal{F}_r}(\Omega;L^2(s,r;U_1))$
satisfies\vspace{-2mm}
\begin{equation}\label{eq:abstract-pathwise-equation}
\dot{\widetilde X}=\mathcal{A}\widetilde X+B\widetilde u+f
\quad\text{in }L^2(\Omega\times(s,r);H).
\end{equation}
Define $u=\Pi_{U_1}\widetilde u\in L^2_{\mathbb F}(s,r;U_1)$. The process $X$ from \cref{eq:conditional-path} has a version in
$L^2_{\mathbb F}(\Omega;C([s,r];H)) \cap L^2_{\mathbb F}(s,r;\mathcal D)$ that agrees with $\Pi_{\mathcal D}\widetilde X$ in $L^2(\Omega\times(s,r);\mathcal D)$. 

Furthermore, there exists a diffusion control $v\in L^2_{\mathbb F}(s,r;\mathcal L_2(\mathfrak K;H))$ such that the state of the control system\vspace{-2mm}
\begin{equation}
\label{eq:abstract-stochastic-equation}
\begin{cases}
d X=(\mathcal{A}X+Bu+f)\,d t+(CX+v)\,d W(t) & t\in (s,r), \\[2pt]
X(s)=\widetilde X(s),  
\end{cases}\vspace{-2mm}
\end{equation}
satisfies $X(r)=\widetilde X(r)$ and\vspace{-2mm}
\begin{equation}
\mathbb{E}\int_s^r\lVert X\rVert_{\mathcal D}^2\,d t
\leq\mathbb{E}\int_s^r\lVert \widetilde X\rVert_{\mathcal D}^2\,d t,
\label{eq:abstract-D-estimate}\vspace{-2mm}
\end{equation}
and the controls satisfy\vspace{-2mm}
\begin{equation} \label{eq:abstract-U-estimate}
\mathbb{E}\int_s^r\lVert u\rVert_{U_1}^2\,d t
\leq\mathbb{E}\int_s^r\lVert \widetilde u\rVert_{U_1}^2\,d t,\vspace{-2mm}
\end{equation}
and\vspace{-1mm}
\begin{equation}\label{eq:abstract-V-estimate}
\mathbb{E}\int_s^r\lVert v\rVert_{\mathcal{L}_2(\mathfrak K;H)}^2\,d t
\leq2(r-s)\mathbb{E}\int_s^r\bigl\lVert \dot{\widetilde X}\bigr\rVert_H^2\,d t
+2\lVert C\rVert_\infty^2\mathbb{E}\int_s^r\lVert \widetilde X\rVert_H^2\,d t.
\end{equation}
Let $K_H\subset H$ and $K_{U_1}\subset U_1$ be fixed nonempty closed convex sets. If $\widetilde X(t)\in K_H$ for all $t\in[s,r]$ almost surely and $\widetilde u(t)\in K_{U_1}$ for almost every $(\omega,t)$, then the same constraints hold for $X$ and $u$.
\end{theorem}
\begin{proof}
\textbf{Step 1. Project the trajectory and identify its $\mathcal{D}$-regularity.}
Set $u=\Pi_{U_1}\widetilde u$.
By \cref{thm:causal-projection}, there exist
$X\in L^2_{\mathbb F}(\Omega;C([s,r];H))$ and
$\Sigma\in L^2_{\mathbb F}(s,r;\mathcal L_2(\mathfrak K;H))$ such that
\begin{equation*}
\begin{cases}
d X=\Pi_H\dot{\widetilde X}\,d t+\Sigma\,d W(t),\\
X(s)=\widetilde X(s),\qquad X(r)=\widetilde X(r).
\end{cases}
\end{equation*}
We next identify the $\mathcal{D}$-valued representative of $X$. Let $\iota:\mathcal{D}\hookrightarrow H$ denote the continuous embedding. Since the predictable projection is the Bochner conditional expectation with respect to the predictable $\sigma$-algebra, it commutes with the bounded linear embedding $\iota$. Hence,
\begin{equation*}
\iota\Pi_{\mathcal D}\widetilde X
=\Pi_H(\iota\widetilde X)=X
\quad\text{in }L^2(\Omega\times(s,r);H).
\end{equation*}
This shows that $X$ agrees with $\Pi_{\mathcal D}\widetilde X$ and thus belongs to $L^2_{\mathbb F}(s,r;\mathcal D)$.

\smallskip

\textbf{Step 2. Identify the projected drift and diffusion terms.} 
To compute the projected drift, we first analyze the action of the operator $\mathcal A$. Since $\mathcal A(t,\omega)\in\mathcal L(\mathcal D;H)$ is bounded for almost every $(\omega,t)$, it admits a Hilbert adjoint $\mathcal A^*(t,\omega)\in\mathcal L(H;\mathcal D)$, defined by 
\begin{equation*}
\langle \mathcal A x, y\rangle_H = \langle x, \mathcal A^* y\rangle_{\mathcal D}, \qquad \forall x\in\mathcal D,\ y\in H.
\end{equation*}
For any test process $\eta\in L^2_{\mathbb F}(s,r;H)$, the strong predictability and uniform boundedness of $\mathcal A^*$ imply that $\mathcal A^*\eta\in L^2_{\mathbb F}(s,r;\mathcal D)$. 

Because $\Pi_{\mathcal D}$ is the orthogonal projection onto $L^2_{\mathbb F}(s,r;\mathcal D)$, the difference $\widetilde X-\Pi_{\mathcal D}\widetilde X$ is orthogonal to any predictable $\mathcal D$-valued process. Therefore,  we obtain
\begin{equation*}
\mathbb E\int_s^r
\bigl\langle\mathcal A\widetilde X-\mathcal A\Pi_{\mathcal D}\widetilde X,
\eta\bigr\rangle_H\,d t
=\mathbb E\int_s^r
\bigl\langle\widetilde X-\Pi_{\mathcal D}\widetilde X,
\mathcal A^*\eta\bigr\rangle_{\mathcal D}\,d t=0.
\end{equation*}
Since $\mathcal A\Pi_{\mathcal D}\widetilde X\in L^2_{\mathbb F}(s,r;H)$, the orthogonality condition above characterizes it as the predictable projection of $\mathcal A\widetilde X$. Hence,\vspace{-2mm}
\begin{equation*}
\Pi_H(\mathcal A\widetilde X)
=\mathcal A\Pi_{\mathcal D}\widetilde X=\mathcal AX.
\end{equation*}

Similarly, we analyze the drift control term $B\widetilde u$. Recall that $B(t,\omega)\in\mathcal L(U_1;H)$ is a strongly predictable and uniformly bounded operator family. Consequently, its Hilbert adjoint $B^*(t,\omega)\in\mathcal L(H;U_1)$ shares the same properties, implying that for any predictable test process $\eta\in L^2_{\mathbb F}(s,r;H)$, the process $B^*\eta$ belongs to $L^2_{\mathbb F}(s,r;U_1)$. Since the predictable projection $\Pi_{U_1}$ is an orthogonal projection, the residual $\widetilde u - \Pi_{U_1}\widetilde u$ is orthogonal to any predictable process. As $B^*\eta$ is predictable, we obtain
\begin{equation*}
\mathbb E\int_s^r \bigl\langle B(t,\omega)\widetilde u - B(t,\omega)\Pi_{U_1}\widetilde u, \eta\bigr\rangle_H \,d t
= \mathbb E\int_s^r \bigl\langle \widetilde u - \Pi_{U_1}\widetilde u, B^*(t,\omega)\eta\bigr\rangle_{U_1} \,d t = 0.
\end{equation*}
Since $B(\cdot)\Pi_{U_1}\widetilde u \in L^2_{\mathbb F}(s,r;H)$, the characterization of the orthogonal projection yields
$
\Pi_H(B\widetilde u) = B\Pi_{U_1}\widetilde u.
$
Recalling that $u = \Pi_{U_1}\widetilde u$, we conclude $\Pi_H(B\widetilde u) = Bu$.

Since the source term $f$ is already predictable, applying $\Pi_H$ to both sides of the pathwise equation \cref{eq:abstract-pathwise-equation} gives
\begin{equation*}
\Pi_H\dot{\widetilde X}
=\Pi_H(\mathcal A\widetilde X)+\Pi_H(B\widetilde u)+f
=\mathcal AX+Bu+f.
\end{equation*}
Thus, by defining the predictable, square-integrable control $v=\Sigma-CX$, the state process $X$ satisfies the  equation \cref{eq:abstract-stochastic-equation}. 

\smallskip

\textbf{Step 3. Prove the bounds and preserve the constraints.}
The contraction properties of $\Pi_{\mathcal D}$ and $\Pi_{U_1}$ give
\cref{eq:abstract-D-estimate,eq:abstract-U-estimate}, respectively. Since $X=\Pi_H\widetilde X$, the contraction property of $\Pi_H$ gives
\begin{equation*}
\mathbb{E}\int_s^r\lVert X\rVert_H^2\,d t
\leq\mathbb{E}\int_s^r\lVert \widetilde X\rVert_H^2\,d t.
\end{equation*}
Combining this inequality with \cref{eq:martingale-defect-estimate} and the definition $v=\Sigma-CX$ proves \cref{eq:abstract-V-estimate}.

Finally, regarding the fixed closed convex constraints $K_H$ and $K_{U_1}$, conditional Jensen's inequality ensures that the projection preserves the control constraint almost everywhere and the state constraint at each fixed time. Since the trajectory $X$ is continuous and $K_H$ is closed, the state constraint extends from a countable dense set of times to all $t\in[s,r]$ almost surely.
\end{proof}
\begin{remark}\label{rem:Vstar-insufficient}
Working on the common domain $\mathcal D = \operatorname{Dom}(\mathcal A_0)$ ensures that the projected drift is $H$-valued and that the martingale integrand belongs to $\mathcal{L}_2(\mathfrak K;H)$. It is natural to ask whether the projection could instead be performed in a larger dual space. For a Hilbert space $V$ continuously and densely embedded in $H$, let $V^*$ denote its dual with $H$ as pivot. Projecting only in $V^*$ would produce a $\mathcal{L}_2(\mathfrak K;V^*)$-valued integrand. However, a diffusion control operator with range in $\mathcal L_2(\mathfrak K;H)$ cannot generally realize such a $V^*$-valued correction. This is why we use the common domain here.
\end{remark}
\begin{remark}\label{rem:convex-causalization-scope}
The fixed closed convex sets may be positive cones, order intervals (i.e., sets bounded above and below by given functions), affine constraints, or control boxes (i.e., pointwise bound constraints on the control).  Nonconvex constraints and sample-dependent
moving sets require a separate argument. The diffusion control inherits no
sign or box constraint from this construction.
\end{remark}

\subsubsection{Realization of the martingale correction and controllability}

The preceding lifting theorem produces a martingale integrand $\Sigma$.
To realize it through the diffusion $CX+Dv$, we must solve
$Dv=\Sigma-CX$. The right-inverse assumption below requires the diffusion
actuator to realize every direction in $\mathcal L_2(\mathfrak K;H)$,
with a predictable solution and a uniform $L^2$ bound.
\begin{corollary}\label{cor:endpoint-transfer}
Assume the hypotheses of \cref{thm:abstract-lifting} on $(s,T)$. Suppose that the diffusion operators 
$C:\Omega\times(s,T)\to\mathcal L(H;\mathcal{L}_2(\mathfrak K;H))$ and
$D:\Omega\times(s,T)\to\mathcal L(U_2;\mathcal{L}_2(\mathfrak K;H))$
are strongly predictable and uniformly bounded, and that
$D$ admits a strongly predictable, uniformly bounded right inverse
$R:\Omega\times(s,T)\to\mathcal L(\mathcal{L}_2(\mathfrak K;H);U_2)$ satisfying
$DR=\operatorname{Id}_{\mathcal{L}_2(\mathfrak K;H)}$.
Then every deterministic controlled trajectory from \cref{thm:abstract-lifting} preserves its endpoint under the stochastic equation with diffusion $CX+Dv$.
Specifically, if $\Sigma$ is the associated martingale integrand, the diffusion control $v=R(\Sigma-CX)$ satisfies
\begin{align}
\lVert v\rVert_{L^2_{\mathbb F}(s,T;U_2)}
&\leq\lVert R\rVert_\infty
\bigl(\lVert \Sigma\rVert_{L^2_{\mathbb F}(s,T;\mathcal{L}_2(\mathfrak K;H))}
+\lVert C\rVert_\infty\lVert X\rVert_{L^2_{\mathbb F}(s,T;H)}\bigr).
\label{eq:master-diffusion-control}
\end{align}
The state and drift control inherit their respective bounds directly from \cref{thm:abstract-lifting}.
\end{corollary}
\begin{proof}
By \cref{thm:abstract-lifting}, there exist an adapted state process $X\in L^2_{\mathbb F}(s,T;H)$, a drift control $u\in L^2_{\mathbb F}(s,T;U_1)$, and a martingale integrand $\Sigma\in L^2_{\mathbb F}(s,T;\mathcal L_2(\mathfrak K;H))$ that satisfy the prescribed endpoints. Define the diffusion control $v:=R(\Sigma-CX)$. Since $C$ and $R$ are strongly predictable, the composition $v = R(\Sigma - CX)$ is strongly predictable. Furthermore, because $R$ and $C$ are uniformly bounded in operator norm and both $X \in L^2_{\mathbb F}(s,T;H)$ and $\Sigma \in L^2_{\mathbb F}(s,T;\mathcal{L}_2(\mathfrak K;H))$, the triangle inequality and the boundedness of $R$ guarantee that $v\in L^2_{\mathbb F}(s,T;U_2)$. Furthermore, using the right-inverse property $DR=\operatorname{Id}_{\mathcal{L}_2(\mathfrak K;H)}$, we obtain
\begin{equation*}
CX + Dv = CX + D(R(\Sigma-CX)) = CX + (\Sigma-CX) = \Sigma.
\end{equation*}
Thus, the state process $X$ satisfies the stochastic evolution equation with the diffusion term $CX+Dv$, and the terminal state is preserved. Finally, applying the triangle inequality and the sub-multiplicative property of the operator norm yields
\begin{align*}
\lVert v\rVert_{L^2_{\mathbb F}(s,T;U_2)}
&\leq \lVert R\rVert_\infty \bigl(\lVert \Sigma\rVert_{L^2_{\mathbb F}(s,T;\mathcal{L}_2(\mathfrak K;H))} + \lVert CX\rVert_{L^2_{\mathbb F}(s,T;\mathcal{L}_2(\mathfrak K;H))}\bigr) \\
&\leq \lVert R\rVert_\infty \bigl(\lVert \Sigma\rVert_{L^2_{\mathbb F}(s,T;\mathcal{L}_2(\mathfrak K;H))} + \lVert C\rVert_\infty\lVert X\rVert_{L^2_{\mathbb F}(s,T;H)}\bigr),
\end{align*}
which proves \cref{eq:master-diffusion-control}. The state and drift control directly inherit the bounds from \cref{thm:abstract-lifting}.
\end{proof}

Since lifting preserves terminal states, exact, null, and approximate
controllability follow from the corresponding deterministic constructions.

\begin{corollary}\label{cor:controllability-transfer}
Assume the coefficient and right-inverse hypotheses of
\cref{cor:endpoint-transfer} on $(s,T)$,
and fix an admissible initial state and source. In each assertion,
assume that the deterministic controlled pairs satisfy all the hypotheses
of \cref{thm:abstract-lifting}.
\begin{enumerate}[label=(\roman*),font=\upshape]
\item If every $\xi\in L^2_{\mathcal F_T}(\Omega;H)$ is the endpoint of
such a trajectory, the stochastic system is exactly controllable from
the given initial state.
\item If such a trajectory ends at zero, the stochastic system is null
controllable from the given initial state.
\item If a sequence of such trajectories has endpoints
$\xi_m\to\xi$ in $L^2_{\mathcal F_T}(\Omega;H)$ as $m\to\infty$, then $\xi$ is  
approximately controllable from the given initial state.  
\end{enumerate}
\end{corollary}

\begin{proof}
Apply \cref{cor:endpoint-transfer} to each controlled trajectory. Since
the lifted terminal state equals the original endpoint, all three
conclusions follow.
\end{proof}

Under the same assumptions, the terminal constraint $X(T)\in\mathcal T$,
where $\mathcal T\subset L^2_{\mathcal{F}_T}(\Omega;H)$, is feasible whenever one
element of $\mathcal T$ is the endpoint of such a controlled trajectory.
The set $\mathcal T$ need not be convex.
Fixed closed convex constraints
on the state and drift control are preserved as in
\cref{thm:abstract-lifting}.

For approximate controllability, assume the hypotheses of
\cref{cor:endpoint-transfer} and suppose the strong-operator
$\mathcal{F}_T$-measurable, uniformly bounded endpoint map
$\Lambda_T:\Omega\to\mathcal{L}(L^2(s,T;U_1);H)$
has dense range almost surely and
$X_T^0\in L^2_{\mathcal{F}_T}(\Omega;H)$. For each target
$\xi\in L^2_{\mathcal{F}_T}(\Omega;H)$, use \cref{eq:common-Tikhonov} with
$L=\Lambda_T$ and $g=\xi-X_T^0$. If the resulting controlled trajectories
satisfy the regularity and integrability hypotheses of
\cref{cor:endpoint-transfer}, their lifts have
terminal values $\xi_m=X_T^0+\Lambda_Tu_m\to\xi$ in
$L^2_{\mathcal F_T}(\Omega;H)$.

The following example shows that lifting a minimum-energy deterministic
control need not minimize the stochastic control energy.

\begin{example}\label{ex:lifted-versus-stochastic-minimum}
Let $W$ be a standard Brownian motion on $[0,1]$. The deterministic
control minimizing $\int_0^1|\widetilde u|^2\,d t$ for
$\dot{\widetilde X}=\widetilde u$, $\widetilde X(0)=0$,
$\widetilde X(1)=W(1)$ is
$\widetilde u(t)=W(1)$, $\widetilde X(t)=tW(1)$.
For the stochastic system $dX=u\,d t+v\,d W$, $X(0)=0$, its lift is
$X(t)=tW(t)$, with $u(t)=W(t)$ and $v(t)=t$. Its energy is
\begin{equation*}
\mathbb{E}\int_0^1\bigl(|u|^2+|v|^2\bigr)\,d t=\frac56.
\end{equation*}
The predictable pair $u^*(t)=\int_0^t(2-s)^{-1}\,d W(s)$,
$v^*(t)=(2-t)^{-1}$ has lower energy:
\begin{equation*}
\mathbb{E}\int_0^1\bigl(|u^*|^2+|v^*|^2\bigr)\,d t
=\int_0^1\frac{\,d s}{2-s}=\log2<\frac56.
\end{equation*}
It reaches the same target, since  
\begin{equation*}
\int_0^1u^*(t)\,d t+\int_0^1v^*(t)\,d W(t)
=\int_0^1\frac{(1-s)+1}{2-s}\,d W(s)=W(1).
\end{equation*}
\end{example}

\section{Null Controllability of Stochastic Parabolic Equations}\label{sec:application-parabolic}

In this and the following two sections, we apply the general framework to three classes of stochastic parabolic equations. We begin with the continuous case, which illustrates how deterministic observability grounds the samplewise deterministic reduction and how the adapted lifting accommodates random coefficients. 
Throughout these three sections, we take $\mathfrak K=\mathbb R$ and let $W$ be a one-dimensional Brownian motion.

\subsection{Formulation of the problem and the controllability result}\label{sec:parabolic-setting}

Let $G\subset\mathbb R^n$ be a bounded domain with a $C^2
$ boundary, and let $G_0\subset G$ be nonempty and open. Set $H=L^2(G)$, $V=H_0^1(G)$, $\mathcal D=H^2(G)\cap H_0^1(G)$, and $U_1=L^2(G_0)$. We equip $V$ and $\mathcal D$ with the $H^1(G)$ and $H^2(G)$ norms, respectively. The operator $B\in\mathcal L(U_1;H)$ is extension by zero, and $D\in\mathcal L(H)$ is the identity operator.

Let $A=(a^{ij})_{1\leq i,j\leq n}$ and $\mathbf b=(b^1,\ldots,b^n)$ be deterministic coefficients on $G$, and let $a_1,a_2$ be real-valued random coefficients on $\Omega\times(0,T)\times G$. These coefficients are assumed to satisfy the following conditions.
\begin{assumption}\label{ass:coefficients}
\begin{enumerate}[label=(\roman*),font=\upshape]
\item The matrix $A$ is symmetric with entries in $C^2(\overline G)$, and $\mathbf b\in W^{1,\infty}(G;\mathbb R^n)$. There exists a constant $\lambda>0$ such that\vspace{-1mm}
\begin{equation*}
\sum_{i,j=1}^n a^{ij}(x)\xi_i\xi_j\geq\lambda|\xi|^2
\quad\text{for all }(x,\xi)\in\overline G\times\mathbb R^n.
\end{equation*}
\item The coefficients $a_1,a_2$ have jointly measurable representatives whose $(\omega,t)$-sections are predictable for almost every $x\in G$. There exist a constant $M>0$ and a set $\Omega_{\rm par}\in\mathcal F_0$ of full probability such that, for every $\omega\in\Omega_{\rm par}$,
\begin{equation*}
\lVert A\rVert_{W^{2,\infty}(G)}
+\lVert \mathbf b\rVert_{W^{1,\infty}(G;\mathbb R^n)}
+\lVert a_1(\omega)\rVert_{L^\infty((0,T)\times G)}
+\lVert a_2(\omega)\rVert_{L^\infty((0,T)\times G)}
\leq M.
\end{equation*}
\end{enumerate}
\end{assumption}
On the $\mathcal F_0$-measurable null set
$\Omega\setminus\Omega_{\rm par}$, set $a_1=a_2=0$.

Define the Dirichlet realization $\mathcal A_0:\mathcal D\to H$ of the elliptic operator by
\begin{equation*}
\mathcal A_0z=\operatorname{div}(A\nabla z)+\mathbf b\cdot\nabla z,
\qquad \operatorname{Dom}(\mathcal A_0)=\mathcal D.
\end{equation*}
It generates an analytic $C_0$-semigroup on $H$, and its graph norm is equivalent to the $H^2(G)$ norm on $\mathcal D$. 
Set $N(t,\omega)z=a_1(t,\omega,\cdot)z$ and $\mathcal A(t,\omega)=\mathcal A_0+N(t,\omega)$. 
The family $N$ is strongly predictable and satisfies $\lVert N\rVert_\infty\leq M$. 
Consequently, $\mathcal A:\mathcal D\to H$ is strongly predictable, and there exists a constant $\mathcal C(M,G)>0$ such that
\begin{equation}\label{eq:parabolic-operator-bound}
\lVert \mathcal A(t,\omega)z\rVert_H
\leq\mathcal C(M,G)\lVert z\rVert_{\mathcal D}
\end{equation}
for almost every $(\omega,t)$. In the variational formulation below, $\mathcal A$ is also identified with its continuous extension as an operator from $V$ to $V^*$.

We consider the controlled stochastic parabolic equation:
\begin{equation}\label{eq:main-system}
\begin{cases}
\displaystyle d y=(\mathcal{A} y+Bu)\,d t+(a_2y+v)\,d W(t),
\quad 0<t<T,\\[2pt]
y(0)=y_0.
\end{cases}
\end{equation}
The homogeneous Dirichlet boundary condition is encoded in $V$ and $\mathcal D$. The drift control $u$ is localized in $G_0$, while the diffusion control $v$ acts on the whole domain $G$.

The following theorem gives null controllability for this system.
For the usual augmented Brownian filtration and domains with $C^4$ boundary,
null controllability in this fixed-generator case follows from \cite[Theorem~2.1]{TangZhang2009}.
Here we obtain it by deterministic observability and adapted lifting.
\begin{theorem}\label{thm:main}
Assume \cref{ass:coefficients}. Then the control system \cref{eq:main-system} is null controllable at every $T>0$. That is, for every $y_0\in L^2_{\mathcal{F}_0}(\Omega;H)$, there exist controls $u\in L^2_{\mathbb F}(0,T;U_1)$ and $v\in L^2_{\mathbb F}(0,T;H)$ such that the corresponding variational solution satisfies $y(T)=0$ almost surely, and
\begin{equation}\label{eq:main-cost}
\mathbb{E}\int_0^T\bigl(\lVert u(t)\rVert_{U_1}^2+\lVert v(t)\rVert_H^2\bigr)\,d t
\leq \mathcal{C}(T)\mathbb{E}\lVert y_0\rVert_H^2.
\end{equation}
The constant $\mathcal{C}(T)$ depends only on $T,G,G_0,\lambda$, and $M$, uniformly over the coefficients in \cref{ass:coefficients}.
\end{theorem}
\begin{remark}\label{rem:random-principal-part}
Theorem~2.1 of \cite{TangZhang2009} allows random, time-dependent principal coefficients. 
The projection step in the proof of \cref{thm:abstract-lifting} also applies to strongly predictable, uniformly bounded families $\widehat{\mathcal A}(t,\omega)\in\mathcal L(\mathcal D;H)$ on a common domain, provided the stochastic equation is well posed and the measurable samplewise controlled trajectories satisfy the stated strong-regularity and integrability bounds. This includes random nonautonomous parabolic operators of Tang--Zhang type when the principal matrix is symmetric and uniformly elliptic, its paths are uniformly bounded in $C^1([0,T];W^{2,\infty}(G))$, the first-order coefficient paths are uniformly bounded in $L^\infty(0,T;W^{1,\infty}(G;\mathbb R^n))$, and the zeroth-order coefficients are bounded and predictable. The uniform bounds are deterministic, and the operator family is strongly predictable on the common domain. Variational well-posedness, deterministic maximal regularity \cite[Corollary~5.2]{Dier2015}, and deterministic observability \cite[Theorem~2.1 and Corollary~2.1]{ImanuvilovYamamoto2003} provide the necessary inputs, and the measurable selection and projection argument then proceed as above.

Although the present technique can handle this more general setting, we restrict ourselves to the fixed principal part to illustrate the core ideas without pursuing technical generality.  
\end{remark}

\subsection{Pathwise regularity and null-control operators}\label{sec:deterministic}

In this subsection, we fix $\omega\in\Omega_{\rm par}$ and treat
the equation with the frozen random potential as a deterministic one.
The energy and maximal regularity estimates are uniform in $\omega$
and over subintervals of $[0,T]$.
The observability and null-control bounds are uniform on intervals
whose lengths are bounded below by a fixed $\delta>0$.
These estimates provide the inputs for the measurable selection below.

Throughout this subsection and in \cref{sec:measurable}, the constants $\mathcal{C}$ in the energy and maximal regularity estimates depend only on $T,G,\lambda$, and $M$. For control estimates and selected trajectories on intervals of length at least a fixed $\delta>0$, the constant $\mathcal{C}(\delta)$ may additionally depend on $G_0$ and $\delta$.

Fix $0\leq s<r\leq T$ and $\omega\in\Omega_{\rm par}$, and suppress $\omega$ from the notation. Consider the  control system
\begin{equation}\label{eq:deterministic-state}
\begin{cases}
\partial_tz-\operatorname{div}(A\nabla z)
=a_1z+\mathbf{b}\cdot\nabla z+Bu+f 
&\text{in }(s,r)\times G,\\[2pt]
z=0 &\text{on }(s,r)\times\partial G,\\[2pt]
z(s)=\xi&\text{in } G.
\end{cases}
\end{equation}
Write $z(\cdot;\xi,u,f)$ for this solution. By \cite[Section~7.1.2, Theorems~2--4]{Evans2010}, for $\xi\in H$, $u\in L^2(s,r;U_1)$, and $f\in L^2(s,r;H)$,
\cref{eq:deterministic-state} has a unique  solution satisfying
\begin{equation}\label{eq:deterministic-energy}
\lVert z\rVert_{C([s,r];H)}^2+\lVert z\rVert_{L^2(s,r;V)}^2
\leq \mathcal C\bigl(\lVert \xi\rVert_H^2+\lVert u\rVert_{L^2(s,r;U_1)}^2
+\lVert f\rVert_{L^2(s,r;H)}^2\bigr).
\end{equation}
Further, by \cite[Theorem~4.1 and Corollary~5.2]{Dier2015}, if $\xi\in V$, then $z\in\mathscr X(s,r)\cap C([s,r];V)$ and
\begin{equation}\label{eq:deterministic-MR}
\lVert z\rVert_{\mathscr X(s,r)}^2+\lVert z\rVert_{C([s,r];V)}^2
\leq\mathcal C\bigl(\lVert \xi\rVert_V^2+\lVert u\rVert_{L^2(s,r;U_1)}^2
+\lVert f\rVert_{L^2(s,r;H)}^2\bigr).
\end{equation}

The next lemma is an immediate corollary of the null controllability result of \eqref{eq:deterministic-state} (see \cite[Theorem~2.1 and Corollary~2.1]{ImanuvilovYamamoto2003}).
\begin{lemma}\label{prop:deterministic-MR}
Assume \cref{ass:coefficients}. Restoring $\omega$, define the bounded endpoint maps $E_\omega:H\to H$ and $\Lambda_\omega:L^2(s,r;U_1)\to H$ by\vspace{-1mm}
\begin{equation*}
E_\omega \xi=z(r;\xi,0,0)
\quad\text{and}\quad
\Lambda_\omega u=z(r;0,u,0).\vspace{-1mm}
\end{equation*}
For every $\xi\in H$, the equation $\Lambda_\omega u=-E_\omega \xi$ admits a unique minimum-norm solution $u_\omega^\dagger(\xi)$. Moreover, the map $\xi\mapsto u_\omega^\dagger(\xi)$ is linear, and for $r-s\geq\delta>0$ there exists a constant $\mathcal C(\delta)>0$ such that
\begin{align}\label{eq:minimum-control-cost}
\lVert u_\omega^\dagger(\xi)\rVert_{L^2(s,r;U_1)}^2
\leq\mathcal C(\delta)\lVert \xi\rVert_H^2.
\end{align}
\end{lemma}

These samplewise estimates are uniform in $\omega$ and provide the boundedness and regularity required by the measurable selection framework. The next subsection makes this connection explicit by verifying the measurability of the associated solution and endpoint operators.

\subsection{Measurable selection of samplewise controlled trajectories}\label{sec:measurable}

The samplewise deterministic estimates established in the previous subsection are uniform in the sample point $\omega$. In this subsection, we combine these estimates with the measurable selection framework from \cref{sec:measurable-common} to obtain measurable, square-integrable controlled trajectories. 

The next lemma verifies the strong-operator $\mathcal{F}_r$-measurability of the solution and endpoint operators required by \cref{thm:random-Tikhonov}.
\begin{lemma}\label{lem:measurable-solution-operators}
Under \cref{ass:coefficients}, the endpoint maps $E_\omega$ and $\Lambda_\omega$
in \cref{prop:deterministic-MR} are strong-operator $\mathcal{F}_r$-measurable and uniformly bounded. For $\xi\in L^2_{\mathcal{F}_r}(\Omega;V)$,
$u\in L^2_{\mathcal{F}_r}(\Omega;L^2(s,r;U_1))$, and
$f\in L^2_{\mathcal{F}_r}(\Omega;L^2(s,r;H))$, the solution of
\cref{eq:deterministic-state} is a strongly $\mathcal{F}_r$-measurable,
square-integrable $\mathscr{X}(s,r)$-valued random variable.
\end{lemma}
\begin{proof}
For each fixed $\omega \in \Omega_{\rm par}$, define the operator
\begin{equation*}
\mathcal T_\omega:\mathscr X(s,r)\longrightarrow L^2(s,r;H)\times V,
\quad\mathcal T_\omega z=(z'-\mathcal A_\omega z,z(s)).
\end{equation*}
Clearly, $\mathcal T_\omega$ is a bounded linear operator, since its components are controlled by the norm of $\mathscr X(s,r)$. 
By \cref{eq:deterministic-MR}, the estimate
\begin{equation*}
\lVert z\rVert_{\mathscr X(s,r)} 
\leq C \bigl( \lVert z'-\mathcal A_\omega z \rVert_{L^2(s,r;H)} + \lVert z(s) \rVert_V \bigr)
\end{equation*}
holds with a constant $C$ independent of $\omega$. This shows that $\mathcal T_\omega$ is injective and its inverse $\mathcal T_\omega^{-1}$ is uniformly bounded. The existence part of the maximal regularity estimate \cref{eq:deterministic-MR} further ensures that $\mathcal T_\omega$ is surjective. Hence, $\mathcal T_\omega$ is a bounded isomorphism with a uniformly bounded inverse.

By \cref{ass:coefficients}, the coefficient $a_1$ is jointly measurable and predictable in $(\omega,t)$. Since the operator $\mathcal A_0$ is deterministic and $z$ is fixed, $a_1z$ is a predictable $H$-valued process on $(s,r)$. Hence $\omega\mapsto \mathcal T_\omega z$ is strongly $\mathcal F_r$-measurable. For measurable data
$F_\omega=(Bu+f,\xi)$, take a countable dense family $\{z_j\}_{j=1}^\infty$ in
$\mathscr X(s,r)$ and choose the first index $j_m(\omega)$ for which
$\lVert\mathcal T_\omega z_j-F_\omega\rVert<1/m$. This index is measurable,
and
\begin{equation*}
\lim_{m\to\infty}\lVert z_{j_m(\omega)}-\mathcal T_\omega^{-1}F_\omega\rVert_{\mathscr X(s,r)}
\leq\lim_{m\to\infty}\lVert\mathcal T_\omega^{-1}\rVert/m=0.
\end{equation*}
Thus the strong trajectory is measurable, and \cref{eq:deterministic-MR}
gives square integrability. Approximating $H$-valued initial data by
$V$-valued data and using \cref{eq:deterministic-energy} gives measurable
weak trajectories in $C([s,r];H)\cap L^2(s,r;V)$. Terminal evaluation proves strong-operator measurability of $E_\omega$ and $\Lambda_\omega$: since the solution $z$ is strongly measurable and the evaluation map $z \mapsto z(r)$ is continuous, the compositions are strongly measurable. Their uniform bounds follow directly from the energy estimate: setting $u=0$ or $\xi=0$ in \cref{eq:deterministic-energy} yields $\lVert E_\omega \rVert \leq \mathcal{C}$ and $\lVert \Lambda_\omega \rVert \leq \mathcal{C}$ for all $\omega$, where $\mathcal{C}$ is independent of $\omega$.
\end{proof}
Fix $r-s\geq\delta>0$. For $\xi\in L^2_{\mathcal{F}_s}(\Omega;V)$,
we apply \cref{thm:random-Tikhonov} with
$L=\Lambda_\omega$, $g=-E_\omega \xi$, and $\mathcal G=\mathcal{F}_r$.
\cref{prop:deterministic-MR} and
\cref{eq:minimum-control-cost} provide the solvability condition and
the square-integrable majorant required by \cref{thm:random-Tikhonov}.
Set $\widetilde u(\omega)=u_\omega^\dagger(\xi(\omega))$ and
$\widetilde y(\omega)=z(\cdot;\xi(\omega),\widetilde u(\omega),0)$, with the coefficients fixed at $\omega$.
The selector and \cref{lem:measurable-solution-operators} give
\begin{equation}\label{eq:selected-control-integrability}
\widetilde u\in L^2_{\mathcal{F}_r}(\Omega;L^2(s,r;U_1)).
\end{equation}
The corresponding trajectory satisfies
\begin{equation}\label{eq:bridge-integrability}
\widetilde y\in L^2_{\mathcal{F}_r}(\Omega;\mathscr{X}(s,r)).
\end{equation}
Its endpoints are $\widetilde y(s)=\xi$ and $\widetilde y(r)=0$, and
\begin{equation}\label{eq:bridge-estimate}
\mathbb{E}\lVert \widetilde u\rVert_{L^2(s,r;U_1)}^2
+\mathbb{E}\lVert \widetilde y\rVert_{\mathscr{X}(s,r)}^2
\leq \mathcal{C}(\delta)\mathbb{E}\lVert \xi\rVert_V^2.
\end{equation}
The resulting pair $(\widetilde y, \widetilde u)$ is a measurable, square-integrable samplewise deterministic trajectory satisfying the exact zero-controllability condition and the required energy bounds. However, because it is selected using the terminal condition and potentially future coefficients, it may fail to be adapted. The adapted lifting established in \Cref{sec:projection} will serve to convert this pair into an admissible stochastic control.

\subsection{Proof of \texorpdfstring{\Cref{thm:main}}{Theorem \ref{thm:main}}}\label{sec:main-proof}

\begin{proof}
\textbf{Step 1. Construct measurable weak paths for $H$-valued initial data.}
Set $\tau=T/2$ and consider \cref{eq:deterministic-state} on $[0,\tau]$ with $u=f=0$ and initial state $y_0 \in L^2_{\mathcal F_0}(\Omega;H)$. Since the strong-regularity result in \cref{lem:measurable-solution-operators} requires $V$-valued initial data, we use the density of $V$ in $H$. Choose a sequence of strongly $\mathcal F_0$-measurable, $V$-valued simple random variables $\{y_{0,m}\}_{m=1}^\infty$ such that $\mathbb E\lVert y_{0,m}-y_0\rVert_H^2\to 0$. For each $m$, let $X_m$ be the corresponding solution. By \cref{lem:measurable-solution-operators}, $X_m$ is strongly measurable in $C([0,\tau];H)$.

By the energy estimate \cref{eq:deterministic-energy}, the solution map is continuous with respect to the initial data:
\begin{equation*}
\mathbb E\sup_{0\leq t\leq \tau}\lVert X_{m_1}(t) - X_{m_2}(t) \rVert_H^2 \leq \mathcal{C} \mathbb E\lVert y_{0,m_1} - y_{0,m_2} \rVert_H^2,\qquad \forall m_1,m_2\in\mathbb N.
\end{equation*}
Thus, $\{X_m\}_{m=1}^\infty$ is Cauchy in the complete space $L^2_{\mathcal F_\tau}(\Omega; C([0,\tau];H))$, and converges to some $X$. The pathwise estimate \cref{eq:deterministic-energy} identifies $X$ as the solution corresponding to $y_0$.

To show that $X$ is adapted, fix any deterministic time $t\in(0,\tau]$. If we restrict the equation to $[0,t]$, the same measurability argument yields a strongly $\mathcal F_t$-measurable $C([0,t];H)$-valued trajectory $X^{[0,t]}$, so its terminal value $X^{[0,t]}(t)$ is $\mathcal F_t$-measurable. By uniqueness, $X(t) = X^{[0,t]}(t)$ almost surely, which proves that $X(t)$ is $\mathcal F_t$-measurable. At $t=0$, $X(0)=y_0$ is $\mathcal F_0$-measurable. Hence $X$ is adapted.

\textbf{Step 2. Regularize the state at the switching time.}
By \cref{eq:deterministic-energy},
\begin{equation}\label{eq:causal-energy}
\mathbb{E}\sup_{0\leq t\leq\tau}\lVert X(t)\rVert_H^2
+\mathbb{E}\int_0^\tau\lVert X(t)\rVert_V^2\,d t
\leq \mathcal{C}(T)\mathbb{E}\lVert y_0\rVert_H^2.
\end{equation}
Choose a deterministic smooth cutoff $\chi$ vanishing near zero and equal to one on $[T/4,\tau]$. Since $(\chi X)_t=\mathcal{A}(\chi X)+\chi'X$ and $(\chi X)(0)=0$, \cref{eq:deterministic-MR} and \cref{lem:measurable-solution-operators} give
$X(\tau)\in L^2_{\mathcal F_\tau}(\Omega;V)$ with
\begin{equation}
\mathbb{E}\lVert X(\tau)\rVert_V^2
\leq \mathcal{C}(T)\mathbb{E}\lVert y_0\rVert_H^2.
\label{eq:prepared-state-bound}
\end{equation}
The controls $(u,v)=(0,-a_2X)$ realize this preparation in \cref{eq:main-system}.

\textbf{Step 3. Select a strong controlled trajectory ending at zero.}
Apply \cref{eq:selected-control-integrability,eq:bridge-integrability,eq:bridge-estimate} on $(\tau,T)$ with $\delta=T/2$ and initial state $X(\tau)$. Together with \cref{eq:prepared-state-bound}, this gives an $\mathcal{F}_T$-measurable control and trajectory with endpoints $\widetilde y(\tau)=X(\tau)$, $\widetilde y(T)=0$, and
\begin{align}\label{eq:pathwise-bridge-cost}
\mathbb{E}\lVert \widetilde y\rVert_{\mathscr{X}(\tau,T)}^2+
\mathbb{E}\lVert \widetilde u\rVert_{L^2(\tau,T;U_1)}^2
\leq \mathcal{C}(T)\mathbb{E}\lVert y_0\rVert_H^2.
\end{align}

\textbf{Step 4. Lift and concatenate.}
By \cref{eq:parabolic-operator-bound}, \cref{thm:abstract-lifting} applies with the generator $\mathcal A_0$ defined above, $N(t,\omega)z=a_1(t,\omega,\cdot)z$, $C(t,\omega)z=a_2(t,\omega,\cdot)z$, $D=\operatorname{Id}_H$, and $f=0$. It supplies $(y_+,u_+,v_+)$ on $(\tau,T)$ with $y_+(\tau)=X(\tau)$ and $y_+(T)=0$.
Concatenating with the preparation yields a variational solution on $[0,T]$, since the states agree at the deterministic switching time. Uniqueness identifies it with the solution of \cref{eq:main-system}.
The preparation estimate \cref{eq:causal-energy} bounds its diffusion control. On $(\tau,T)$, combine \cref{eq:pathwise-bridge-cost} with \cref{eq:abstract-U-estimate,eq:abstract-V-estimate}. Adding the two interval costs gives \cref{eq:main-cost}.
\end{proof}

\section{Mesh-independent null controllability of semi-discrete stochastic parabolic equations}\label{sec:application-semidiscrete}

In this section, we apply the general framework to the spatial semi-discretization of stochastic parabolic equations. A key requirement for numerical approximation is that the controllability property, along with its control cost, is preserved uniformly as the mesh size $h$ tends to zero. We prove null controllability for every sufficiently fine mesh with a cost independent of $h$. Notably, the underlying deterministic system may contain a bounded diagonal potential that is only measurable in time.

Our approach combines the samplewise deterministic reduction and the adapted lifting established in \Cref{sec:framework}. On the deterministic side, a relaxed control estimate \cite{Boyer2014} first yields an exponentially small terminal remainder, which is then eliminated by a short-time deterministic correction. This provides an exact null control with a mesh-independent cost. On the stochastic side, the standard maximal regularity estimate would involve the mesh-dependent operator norm of the discrete Laplacian. To bypass this obstacle, we estimate the effective martingale integrand using the discrete $A_h^{1/2}$-energy, ensuring that the lifted stochastic control cost remains uniformly bounded with respect to $h$.

For autonomous one-dimensional equations with continuous potentials, uniform null controllability was proved in \cite[Theorem~5.3, p.~902]{Allonsius2018}. Lecaros, P\'erez, and Prado \cite{Lecaros2025} established relaxed controllability in arbitrary dimensions but also exposed spatial discretizations with unobservable modes for particular control geometries. In contrast, our result handles nonautonomous, time-dependent diagonal potentials and establishes exact stochastic null controllability with a strictly mesh-independent cost.

\subsection{The semi-discrete stochastic system}

Let $T>0$ and let $\mathcal O\Subset(0,1)$ be a nonempty open interval.
For $N\geq1$, set $h=(N+1)^{-1}$ and $x_j=jh$ for
$0\leq j\leq N+1$.
For a discrete state vector $y=(y_1,\ldots,y_N)$, set $y_0=y_{N+1}=0$ and define
\begin{equation*}
(A_hy)_j=\frac{2y_j-y_{j-1}-y_{j+1}}{h^2},
\qquad 1\leq j\leq N.
\end{equation*}
Set $\mathcal O_h=\{j\in\{1,\ldots,N\}:x_j\in\mathcal O\}$.
Equip $H_h=\mathbb R^N$ and $\mathcal U_h=\mathbb R^{\mathcal O_h}$ with
the weighted inner product\vspace{-1mm}
\begin{equation*}
\langle X,y\rangle_h=h\sum_{j=1}^NX_jy_j\vspace{-1mm}
\end{equation*}
and control norm\vspace{-1mm}
\begin{equation*}
\lVert u\rVert_{\mathcal U_h}^2=h\sum_{j\in\mathcal O_h}\lvert u_j\rvert^2.\vspace{-1mm}
\end{equation*}
Write $\lVert y\rVert_h^2=\langle y,y\rangle_h$.
The matrix $A_h$ is positive self-adjoint; $A_h^{1/2}$ denotes its positive
square root. The operator $B_h:\mathcal U_h\to H_h$ is the canonical embedding that extends a control vector by zero to the uncontrolled nodes. In particular,
$\lVert B_h\rVert\leq1$.
Throughout this section, mesh sizes are restricted to $h=(N+1)^{-1}$.

We study the semi-discrete stochastic control system
\begin{equation}\label{sd:eq:stochastic-system}
\begin{cases}
d y_h+A_hy_h\,d t
=\bigl(Q_h(t)y_h+B_hu_h\bigr)\,d t
+\bigl(C_h(t)y_h+v_h\bigr)\,d W(t), & t\in [0,T],\\[2pt]
y_h(0)=y_h^0.
\end{cases}
\end{equation}
Here $y_h$ is the grid state with initial value $y_h^0$, and
$u_h\in L^2_{\mathbb F}(0,T;\mathcal U_h)$ and
$v_h\in L^2_{\mathbb F}(0,T;H_h)$ are the controls.
For fixed $M_1,M_2>0$,
assume that $Q_h(t)=\operatorname{diag}(q_{1,h}(t),\ldots,q_{N,h}(t))$
is a predictable diagonal process with
\begin{equation}\label{sd:eq:Q-bound}
\operatorname*{ess\,sup}_{(\omega,t)\in\Omega\times(0,T)}
\max_{1\leq j\leq N}\lvert q_{j,h}(\omega,t)\rvert\leq M_1,
\end{equation}
and $C_h$ is predictable with
\begin{equation}\label{sd:eq:C-bound}
\operatorname*{ess\,sup}_{(\omega,t)\in\Omega\times(0,T)}
\lVert C_h(\omega,t)\rVert_{\mathcal L(H_h)}\leq M_2.
\end{equation}
Both bounds are independent of $h$. Since the mesh family is countable,
there is $\Omega_{\rm sd}\in\mathcal F_0$ of probability one on which these
bounds hold for almost every $t$ and every admissible $h$.
Set $Q_h=C_h=0$ on $\Omega\setminus\Omega_{\rm sd}$.

For each fixed mesh, this system fits the semigroup framework with
$H=H_h$, $\mathcal A_0=-A_h$, $S_0(t)=e^{-tA_h}$,
$N(t,\omega)=Q_h(t,\omega)$, $B=B_h$, $C=C_h$,
and $D=\operatorname{Id}_{H_h}$.

The coefficients are defined directly on the grid. If they are obtained
by sampling continuum coefficients, jointly measurable pointwise
representatives must be specified, since an $L^\infty(G)$ equivalence
class does not determine nodal values.

\subsection{Samplewise deterministic relaxed controls}\label{sec:deterministic-relaxed}

Following the samplewise deterministic reduction, we first fix the sample point $\omega$ and consider the corresponding deterministic equation with the diagonal potential $Q(t)$. Our goal is to establish a deterministic control estimate that is uniform in the mesh size $h$ and the potential $Q$. We begin with a relaxed control estimate (yielding an exponentially small terminal remainder), which will be converted into an exact null control in the next subsection.

For a fixed mesh, we consider the deterministic system
\begin{equation}\label{sd:eq:deterministic-system}
\begin{cases}
y_h'+A_hy_h=Q(t)y_h+B_hu_h, & t\in [0,T],\\[2pt]
y_h(0)=y_h^0.
\end{cases}
\end{equation}
We identify $Q(t)=\operatorname{diag}(q_1(t),\ldots,q_N(t))$ with $q=(q_1,\ldots,q_N)$ and, for $M>0$, define the bounded admissible set
\begin{equation*}
\mathcal Q_h(M):=
\bigl\{
q\in L^1(0,T;\mathbb{R}^N):
\lvert q_j(t)\rvert\leq M
\text{ for a.e. $t$ and every $j$}
\bigr\}.
\end{equation*}
The set $\mathcal Q_h(M)$ is closed in $L^1(0,T;\mathbb R^N)$, since its componentwise bound passes to an almost-everywhere convergent subsequence. We equip it with the relative $L^1$ topology and its Borel $\sigma$-field.

The following proposition recalls the relaxed control estimate of \cite[Propositions~4.1 and 4.3]{Boyer2014}, together with its energy and measurability consequences. The measurability of the control with respect to $Q$ is crucial for the subsequent samplewise measurable selection.

\begin{proposition}
\label{sd:prop:relaxed-control}
There are $\mathcal{C},c,h_0>0$, depending only on $T,\mathcal O$, and $M$, such that for every $0<h\leq h_0$, every $Q\in\mathcal Q_h(M)$, and every $y_h^0\in H_h$, \cref{sd:eq:deterministic-system} has a control $u_h\in L^2(0,T;\mathcal U_h)$ satisfying the relaxed control cost
\begin{equation}
\int_0^T\lVert u_h(t)\rVert_{\mathcal U_h}^2\,d t
\leq \mathcal{C}\lVert y_h^0\rVert_h^2,
\label{sd:eq:relaxed-cost}
\end{equation}
and the exponentially small terminal remainder
\begin{equation}
\lVert y_h(T)\rVert_h
\leq \mathcal{C}e^{-c/h}\lVert y_h^0\rVert_h.
\label{sd:eq:relaxed-terminal}
\end{equation}
The corresponding state also satisfies the mesh-uniform energy estimate
\begin{equation}
\sup_{0\leq t\leq T}\lVert y_h(t)\rVert_h^2
+\int_0^T\lVert A_h^{1/2}y_h(t)\rVert_h^2\,d t
\leq \mathcal{C}\lVert y_h^0\rVert_h^2.
\label{sd:eq:relaxed-energy}
\end{equation}
For each fixed mesh, the control and state can be chosen linearly in $y_h^0$ and Borel measurably in $Q\in\mathcal Q_h(M)$, with values in $L^2(0,T;\mathcal U_h)$ and $C([0,T];H_h)$, respectively.
\end{proposition}
\begin{proof}
Apply \cite[Propositions~4.1 and~4.3, pp.~1057, 1059]{Boyer2014} with potential $a=-Q$. For fixed $T$ and $\lVert Q\rVert_\infty\leq M$, the relaxed observability estimate and penalized Hilbert Uniqueness Method (HUM) yield the relaxed control cost \cref{sd:eq:relaxed-cost} and the exponentially small terminal remainder \cref{sd:eq:relaxed-terminal} uniformly in $h$. The mesh-uniform energy estimate \cref{sd:eq:relaxed-energy} then follows directly from the energy identity for \cref{sd:eq:deterministic-system} and the control cost bound. 

Regarding the measurability and linearity, the penalized HUM control is constructed via the inverse of a regularized Gramian. For each fixed $h$, the evolution of \cref{sd:eq:deterministic-system} depends continuously on $Q$ in the relative $L^1$ topology. Consequently, both the Gramian and its regularized inverse are continuous in $Q$. This implies that the selected control and its trajectory are Borel measurable functions of $Q\in\mathcal Q_h(M)$. Moreover, since the HUM control is obtained by applying the inverse Gramian to the initial data, it is linear with respect to $y_h^0$, which completes the proof.
\end{proof}

As noted in \cref{sd:eq:relaxed-terminal}, the relaxed control leaves an exponentially small remainder at the terminal time $T$. In the next subsection, we will use the short-time estimates in \Cref{app:semidiscrete-technical} to remove this remainder via a deterministic correction, thereby obtaining an exact deterministic null control before applying the stochastic lifting.

\subsection{Exact samplewise deterministic null controllability}\label{sec:deterministic-exact}

The short-time estimate from \Cref{app:semidiscrete-technical} now allows us to remove the exponentially small remainder left by the relaxed control while retaining a mesh-uniform control cost. This yields exact samplewise deterministic null controllability.

To apply the short-time correction, we choose a suitable observation geometry. Choose $x_*\in\mathcal O$ and $d_*>0$ such that $[x_*-3d_*,x_*+3d_*]\subset\mathcal O$. For all sufficiently small $h$, choose $m=m(h)$ so that the adjacent nodes $x_m,x_{m+1}$ belong to $[x_*-d_*,x_*+d_*]$. The left tail $\{1,\ldots,m\}$ and the reversed right tail $\{m+1,\ldots,N\}$ have block lengths $(m+1)h$ and $(N-m+1)h$ in a fixed interval $[L_-,L_+]\subset(0,\infty)$, with bounds depending only on $\mathcal O$.

For a block of $n\geq1$ nodes, we use the standard second-difference matrix with homogeneous Dirichlet boundary conditions
$$A_{h,n}=h^{-2}\operatorname{tridiag}(-1,2,-1)\in\mathbb R^{n\times n}$$ acting on vectors equipped with the weighted norm $\lVert z\rVert_{h,n}^2=h\sum\limits_{j=1}^n\lvert z_j\rvert^2$. Observing these two adjacent nodes controls both grid blocks. The following estimate quantifies the resulting short-time cost uniformly over bounded time-dependent diagonal potentials.
\begin{proposition}
\label{sd:prop:short-observability}
For fixed $M>0$, there exist constants $\mathcal{C},h_*>0$, depending only on $\mathcal O$ and $M$, such that for any $0<h\leq h_*$ and $0<\delta\leq1$, the following holds. Let $Q_h(t)=\operatorname{diag}(q_{1,h}(t),\ldots,q_{N,h}(t))$ on $(0,\delta)$, where $\max_j\lVert q_{j,h}\rVert_{L^\infty(0,\delta)}\leq M$. Every solution of
\begin{equation*}
z_h'+A_hz_h=Q_h(t)z_h
\qquad\text{on }(0,\delta)
\end{equation*}
satisfies the observability estimate
\begin{equation}\label{sd:eq:short-observability}
\lVert z_h(\delta)\rVert_h^2
\leq \mathcal{C}h^{-10}
\exp\biggl(\frac{\mathcal{C}}{\delta}+\frac{\mathcal{C}M\delta}{h}\biggr)
\int_0^\delta\lVert B_h^*z_h(t)\rVert_{\mathcal U_h}^2\,d t.
\end{equation}
By duality, the system is null controllable on $(0,\delta)$ with a squared control cost bounded by the multiplicative factor in \cref{sd:eq:short-observability} times the squared initial norm. Both the estimate and the controllability property remain valid on any translated interval of length $\delta$.
\end{proposition}
\begin{proof}
\textbf{Step 1. Estimate the left and right blocks.}
Write $e_j$ for the standard coordinate vectors in each grid block and
$z_h=(z_1,\ldots,z_N)$. Set $z_L=(z_1,\ldots,z_m)$ and
$Q_L(t)=\operatorname{diag}(q_{1,h}(t),\ldots,q_{m,h}(t))$.
The left block satisfies
\begin{equation*}
	z_L'+A_{h,m}z_L=Q_L(t)z_L+h^{-2}e_mz_{m+1}.
\end{equation*}
Write $z_L=\widehat z_L+w$, where $\widehat z_L$ solves the homogeneous equation on the left block with initial value $z_L(0)$, while $w(0)=0$ and $w$ is driven by the forcing term $h^{-2}e_mz_{m+1}$. By
\cref{sd:prop:robust-endpoint},
\begin{equation}\label{sd:eq:left-homogeneous-observation}
	\lVert \widehat z_L(\delta)\rVert_{h,m}^2
	\leq \mathcal{C}h^{-5}\exp\biggl(\frac{\mathcal{C}}{\delta}+4mM\delta\biggr)
	\int_0^\delta h\lvert \widehat z_{L,m}(t)\rvert^2\,d t.
\end{equation}
The evolution generated by $-A_{h,m}+Q_L(t)$ has norm at most
$e^{M(t-s)}$ in the weighted norm $\lVert\cdot\rVert_{h,m}$.
Duhamel's formula therefore gives, for any $t\in(0,\delta]$,
\begin{equation}\label{sd:eq:left-forcing-terminal}
	\lVert w(t)\rVert_{h,m}^2
	\leq t e^{2Mt}h^{-4}
	\int_0^t h\lvert z_{m+1}(s)\rvert^2\,d s.
\end{equation}
Integrating this inequality over $(0,\delta)$ yields
\begin{equation}\label{sd:eq:left-forcing-bounds}
	\int_0^\delta\lVert w(t)\rVert_{h,m}^2\,d t
	\leq \delta^2 e^{2M\delta}h^{-4}
	\int_0^\delta h\lvert z_{m+1}(s)\rvert^2\,d s.
\end{equation}
Using $h\lvert w_m\rvert^2\leq\lVert w\rVert_{h,m}^2$ and $\widehat z_{L,m}=z_m-w_m$, the homogeneous observation bound \cref{sd:eq:left-homogeneous-observation} yields
\begin{equation*}
	\lVert \widehat z_L(\delta)\rVert_{h,m}^2
	 \leq \mathcal{C}h^{-5}\exp\biggl(\frac{\mathcal{C}}{\delta}+4mM\delta\biggr)
	\int_0^\delta \bigl(2h\lvert z_m(t)\rvert^2+2\lVert w(t)\rVert_{h,m}^2\bigr)\,d t.
\end{equation*}
Substituting the time-integrated forcing bound \cref{sd:eq:left-forcing-bounds} into the above inequality and applying the terminal bound \cref{sd:eq:left-forcing-terminal} to $\lVert w(\delta)\rVert_{h,m}^2$, the triangle inequality $\lVert z_L(\delta)\rVert_{h,m}^2 \leq 2\lVert \widehat z_L(\delta)\rVert_{h,m}^2 + 2\lVert w(\delta)\rVert_{h,m}^2$ gives
\begin{align}
	\lVert z_L(\delta)\rVert_{h,m}^2
	&\leq \mathcal{C}h^{-5}\exp\biggl(\frac{\mathcal{C}}{\delta}+4mM\delta\biggr)
	\int_0^\delta h\lvert z_m(t)\rvert^2\,d t + \mathcal{C}h^{-9}\exp\biggl(\frac{\mathcal{C}}{\delta}+4mM\delta\biggr)
	\int_0^\delta h\lvert z_{m+1}(t)\rvert^2\,d t  \notag\\
	&\quad + \mathcal{C}h^{-4}e^{2M\delta}
	\int_0^\delta h\lvert z_{m+1}(t)\rvert^2\,d t. \label{sd:eq:term3}
\end{align}
Since $m \leq L_+/h$, the exponent $4mM\delta$ is absorbed by $\mathcal{C}M\delta/h$. Enlarging the polynomial factors in \cref{sd:eq:term3} to $h^{-9}$ and factoring out the common coefficient, we combine the integrands into a single integral to obtain
\begin{equation}\label{sd:eq:left-tail-observation}
	\lVert z_L(\delta)\rVert_{h,m}^2
	\leq \mathcal{C}h^{-9}
	\exp\biggl(\frac{\mathcal{C}}{\delta}+\frac{\mathcal{C}M\delta}{h}\biggr)
	\int_0^\delta h
	\bigl(\lvert z_m(t)\rvert^2+\lvert z_{m+1}(t)\rvert^2\bigr)\,d t.
\end{equation}
Due to the symmetry of the grid, the reversed right block (of length $N-m$) satisfies the same type of block equation as the left block, with the forcing term originating from $z_m$ and the observed endpoint being $z_{m+1}$. Applying the left-tail estimate \cref{sd:eq:left-tail-observation} to this reversed block yields an identical bound for $\lVert (z_{m+1},\ldots,z_N)(\delta)\rVert_{h,N-m}^2$, which involves the same observation integrals over $z_m$ and $z_{m+1}$. 

Adding the estimates for the left and right blocks, and using the fact that both observed nodes lie in the control region $\mathcal O_h$ (so that $\int_0^\delta h(|z_m|^2+|z_{m+1}|^2)\,dt \leq \int_0^\delta \lVert B_h^*z_h(t)\rVert_{\mathcal U_h}^2\,dt$), we bound the total energy $\lVert z_h(\delta)\rVert_h^2$. Finally, $h^{-9}\leq h^{-10}$ for $0<h\leq1$, which gives the desired short-time observability estimate
\cref{sd:eq:short-observability}.

\smallskip

\textbf{Step 2. Deduce the control estimate by duality.}
For the control system \cref{sd:eq:deterministic-system}, the adjoint
satisfies $-p'+A_hp=Q_h(t)p$. Thus $\widehat p(r)=p(\delta-r)$ satisfies
$\widehat p'+A_h\widehat p=Q_h(\delta-r)\widehat p$, with
$\widehat p(0)=p(\delta)$ and $\widehat p(\delta)=p(0)$.
The reversed potential has the same diagonal bound, and time reversal
preserves the observation integral. Applying the estimate to $\widehat p$, the HUM for this finite-dimensional system yields the desired control assertion. 

Finally, because the operator $A_h$ is time-independent and the potential bound $M$ holds uniformly for all times, translating the time interval from $(0,\delta)$ to any arbitrary interval of length $\delta$ leaves the constants and exponential factors unchanged. Hence, the observability estimate and the resulting null controllability property remain valid on any such translated interval.
\end{proof}
The short-time estimate removes the exponentially small remainder left by
\cref{sd:prop:relaxed-control} while retaining a mesh-uniform control cost.
\begin{theorem}
\label{sd:thm:deterministic-exact}
For every $M>0$, there exist constants $\mathcal{C},h_0>0$, depending only on $T,\mathcal{O},M$, such that for every $0<h\leq h_0$, every $Q_h\in\mathcal{Q}_h(M)$, and every $y_h^0\in H_h$, the system \cref{sd:eq:deterministic-system} is exactly null controllable. That is, there exists a control $u_h\in L^2(0,T;\mathcal{U}_h)$ whose corresponding state satisfies $y_h(T)=0$. 

Moreover, the control cost and the trajectory satisfy the mesh-uniform bounds
\begin{equation}
\int_0^T\lVert u_h(t)\rVert_{\mathcal U_h}^2\,d t
\leq \mathcal{C}\lVert y_h^0\rVert_h^2, \label{sd:eq:det-exact-cost}
\end{equation}
and
\begin{equation}
\sup_{0\leq t\leq T}\lVert y_h(t)\rVert_h^2
+\int_0^T\lVert A_h^{1/2}y_h(t)\rVert_h^2\,d t
\leq \mathcal{C}\lVert y_h^0\rVert_h^2. \label{sd:eq:det-exact-energy}
\end{equation}
\end{theorem}
\begin{proof}
Apply \cref{sd:prop:relaxed-control} on $(0,T/2)$ to obtain a control with squared cost at most $\mathcal{C}\lVert y_h^0\rVert_h^2$ and $\lVert y_h(T/2)\rVert_h\leq \mathcal{C}e^{-\gamma/h}\lVert y_h^0\rVert_h$ for some $\gamma>0$ independent of $h$. After decreasing $h_0$ if necessary, set $\delta_h=\sqrt h<\min\{T/2,1\}$ and $T_h=T-\delta_h$. Using zero control on $(T/2,T_h)$ and applying the energy inequality, we obtain the exponentially small terminal remainder $\lVert y_h(T_h)\rVert_h\leq \mathcal{C}(T,M)e^{-\gamma/h}\lVert y_h^0\rVert_h$.

Next, we remove this remainder by applying the null-control assertion of \cref{sd:prop:short-observability} on $(T_h,T)$. Denoting the resulting correction by $u_h^{\rm c}$, it gives $y_h(T)=0$ and
\begin{equation*}
\int_{T_h}^T\lVert u_h^{\rm c}\rVert_{\mathcal U_h}^2\,d t
\leq \mathcal{C}h^{-10}\exp\biggl(-\frac{2\gamma}{h}
+\frac{(1+M)\mathcal{C}}{\sqrt h}\biggr)\lVert y_h^0\rVert_h^2
\leq \mathcal{C}\lVert y_h^0\rVert_h^2.
\end{equation*}
The last inequality holds because $h^{-1}$ dominates $h^{-1/2}$ and $|\log h|$ as $h\to 0$. Finally, concatenating the two controls on $(0,T)$ yields a valid control reaching zero at time $T$. Adding the two control costs proves \cref{sd:eq:det-exact-cost}, and the energy estimate for the concatenated equation on $(0,T)$ gives \cref{sd:eq:det-exact-energy}.
\end{proof}
\begin{remark}\label{sd:rem:exponent-balance}
For any $\delta_h=h^\alpha$, $0<\alpha<1$, both
$\delta_h^{-1}$ and $\delta_h/h$ are $o(h^{-1})$ as $h\to0$.
$\alpha=1/2$ balances these two contributions to the correction cost.
\end{remark}

\paragraph{Measurable selection of the minimum-norm null control.}
For each $h\leq h_0$ from \cref{sd:thm:deterministic-exact}, a Gramian
formula selects the minimum-norm control measurably. Let $S_Q(t,s)$ be
the evolution generated by $-A_h+Q(t)$ and set
\begin{equation*}
\mathcal G_Q=\int_0^T S_Q(T,s)B_hB_h^*S_Q(T,s)^*\,d s.
\end{equation*}
Null controllability and invertibility of $S_Q(T,0)$ make $\mathcal G_Q$
positive definite. The minimum-norm null control is
\begin{equation}\label{sd:eq:minimal-control}
(\mathcal K_h(Q)y^0)(s)
=-B_h^*S_Q(T,s)^*\mathcal G_Q^{-1}S_Q(T,0)y^0.
\end{equation}
For fixed $h$, continuous dependence of $S_Q$ on $Q$ in the relative
$L^1$ topology and continuity of matrix inversion imply that
$\mathcal K_h:\mathcal Q_h(M)\to
\mathcal L(H_h;L^2(0,T;\mathcal U_h))$ is continuous in operator norm,
hence Borel measurable. Since $\mathcal K_h(Q)y^0$ is defined as the minimum-norm null control, its control cost is bounded by that of any other null control. Combining this minimality with the existence of a bounded-cost control established in \cref{sd:eq:det-exact-cost}, we obtain
\begin{equation}\label{sd:eq:Kh-bound}
\lVert\mathcal K_h(Q)\rVert_{\mathcal L(H_h;L^2(0,T;\mathcal U_h))}
\leq\mathcal C
\end{equation}
uniformly over $Q\in\mathcal Q_h(M)$ and $h\leq h_0$, with $\mathcal C$
depending only on $T,\mathcal O,M$.
Since $Q_h$ is a predictable process, it can be viewed as a strongly $\mathcal F_T$-measurable random variable taking values in $L^1(0,T;\mathbb{R}^N)$. Consequently, after the modification outside $\Omega_{\rm sd}$, $Q_h$ belongs to $L^0(\Omega; \mathcal Q_h(M_1))$. 
This establishes that the selected samplewise deterministic pair $(\widetilde y_h, \widetilde u_h=\mathcal K_h(Q_h)y_h^0)$ is exactly null-controlling and satisfies uniform deterministic bounds, setting the stage for the stochastic lifting in the next subsection.

\subsection{Mesh-independent stochastic null controllability}\label{sec:stochastic-exact}

We now lift the exact samplewise deterministic null controls established in the previous subsection to the semi-discrete stochastic system \eqref{sd:eq:stochastic-system}. The main challenge here is that the standard abstract lifting theorem (\cref{thm:abstract-lifting}) uses the maximal regularity norm of the trajectory, which involves the operator norm of $A_h$. Since $\lVert A_h\rVert\asymp h^{-2}$, a direct application would give a mesh-dependent upper bound for the diffusion control cost as $h \to 0$.

To obtain a control cost independent of the mesh size, we instead estimate the effective martingale integrand through the discrete $A_h^{1/2}$-energy. This relies on the fact that the projection onto predictable processes commutes with the deterministic operator $A_h^{1/2}$, allowing us to bound the lifted dissipation by the deterministic one. We formulate this lifting procedure in the following lemma.
\begin{lemma}
\label{sd:lem:mesh-lifting}
Let $y_h^0\in L^2_{\mathcal F_0}(\Omega;H_h)$ and $\widetilde u_h\in L^2_{\mathcal F_T}(\Omega;L^2(0,T;\mathcal U_h))$. Let $\widetilde y_h\in L^2_{\mathcal F_T}(\Omega;C([0,T];H_h))$ be the corresponding solution of \cref{sd:eq:deterministic-system} with $Q=Q_h(\omega)$ almost surely. Set $R_h=\widetilde y_h(T)$. Assume there exists a constant $\mathcal C_0$, independent of $h$, such that
\begin{equation}\label{sd:eq:bridge-uniform-bound}
\mathbb{E}\int_0^T\lVert \widetilde u_h\rVert_{\mathcal U_h}^2\,d t
+\mathbb{E}\sup_{0\leq t\leq T}\lVert \widetilde y_h(t)\rVert_h^2
+\mathbb{E}\int_0^T\lVert A_h^{1/2}\widetilde y_h\rVert_h^2\,d t
\leq \mathcal C_0\mathbb{E}\lVert y_h^0\rVert_h^2.
\end{equation}
Then, the stochastic system \cref{sd:eq:stochastic-system} admits predictable controls $(u_h,v_h)$ such that the corresponding state satisfies $y_h(T)=R_h$ almost surely. Moreover, these controls satisfy the mesh-uniform bound
\begin{equation}
\mathbb{E}\int_0^T
\bigl(\lVert u_h\rVert_{\mathcal U_h}^2+\lVert v_h\rVert_h^2\bigr)\,d t
\leq \mathcal{C}\mathbb{E}\lVert y_h^0\rVert_h^2,
\label{sd:eq:lifted-control-bound}
\end{equation}
where the constant $\mathcal{C}$ depends only on $T,M_1,M_2,\mathcal C_0$ and is independent of $h$.
\end{lemma}
\begin{proof}
For each fixed $h$, the deterministic control system \cref{sd:eq:deterministic-system} yields $\widetilde y_h\in L^2_{\mathcal{F}_T}(\Omega;H^1(0,T;H_h))$.
By \cref{thm:causal-projection},
$y_h=\Pi_{H_h}\widetilde y_h$ has a continuous adapted version and,
with $u_h=\Pi_{\mathcal U_h}\widetilde u_h$ and an effective martingale
integrand $G_h$, satisfies
\begin{equation*}
\begin{cases}
d y_h+A_hy_h\,d t=(Q_hy_h+B_hu_h)\,d t+G_h\,d W(t),\qquad t\in (0,T),\\[2pt]
y_h(0)=y_h^0,\quad y_h(T)=R_h.
\end{cases}
\end{equation*}
Since $u_h=\Pi_{\mathcal U_h}\widetilde u_h$, the contraction property of the orthogonal projection implies that
\begin{equation*}
\mathbb{E}\int_0^T\lVert u_h\rVert_{\mathcal U_h}^2\,d t
\leq \mathbb{E}\int_0^T\lVert \widetilde u_h\rVert_{\mathcal U_h}^2\,d t.
\end{equation*} 
To estimate $G_h$ uniformly in $h$,
we apply It\^o's formula to $\lVert y_h\rVert_h^2$:
\begin{align}
\lVert R_h\rVert_h^2-\lVert y_h^0\rVert_h^2
&=-2\int_0^T\lVert A_h^{1/2}y_h\rVert_h^2\,d t
+2\int_0^T\langle Q_hy_h+B_hu_h,y_h\rangle_h\,d t \notag
\\
&\quad+\int_0^T\lVert G_h\rVert_h^2\,d t
+2\int_0^T\langle y_h,G_h\rangle_h\,d W(t).\label{11.22-eq1}
\end{align}
For each fixed $h$, the established $L^2(\Omega;C([0,T];H_h))$ bound for
$y_h$ and $G_h\in L^2_{\mathbb F}(0,T;H_h)$ imply
$\langle y_h,G_h\rangle_h\in L^1(\Omega;L^2(0,T))$, so the stochastic
integral has zero expectation by \cite[Theorem~2.130, item~3]{Lue2021a}.
Taking expectations of \eqref{11.22-eq1} gives
\begin{equation}\label{sd:eq:effective-diffusion-energy}
\mathbb{E}\int_0^T\lVert G_h\rVert_h^2\,d t
\leq\mathbb{E}\lVert R_h\rVert_h^2
+2\mathbb{E}\int_0^T\lVert A_h^{1/2}y_h\rVert_h^2\,d t
+\mathcal{C}\mathbb{E}\int_0^T
\bigl(\lVert y_h\rVert_h^2+\lVert u_h\rVert_{\mathcal U_h}^2\bigr)\,d t.
\end{equation}
Recall that $y_h = \Pi_{H_h}\widetilde y_h$ is the predictable projection of $\widetilde y_h$. Since $A_h^{1/2}$ is a deterministic operator, it commutes with the projection, which implies that $A_h^{1/2}y_h = \Pi_{H_h}(A_h^{1/2}\widetilde y_h)$. By the contraction property of the orthogonal projection, we obtain
\begin{equation}\label{sd:eq:projected-dissipation}
\mathbb{E}\int_0^T\lVert A_h^{1/2}y_h\rVert_h^2\,d t
\leq\mathbb{E}\int_0^T\lVert A_h^{1/2}\widetilde y_h\rVert_h^2\,d t.
\end{equation}
The terminal value $R_h=\widetilde y_h(T)$ is controlled by the
trajectory supremum in \cref{sd:eq:bridge-uniform-bound}. Combining
\cref{sd:eq:bridge-uniform-bound,sd:eq:effective-diffusion-energy,sd:eq:projected-dissipation}
gives an $h$-independent bound for $G_h$.
Setting $v_h=G_h-C_hy_h$ and using \cref{sd:eq:C-bound} proves
\cref{sd:eq:lifted-control-bound}.
\end{proof}

We now apply the preceding lemma to the measurable minimum-norm null controls
in \cref{sd:eq:minimal-control} constructed in the deterministic setting. Their deterministic cost and energy bounds guarantee that the stochastic null controllability is preserved uniformly as $h \to 0$.
\begin{theorem}
\label{sd:thm:stochastic-exact}
Assume \cref{sd:eq:Q-bound,sd:eq:C-bound}. There exist constants $\mathcal{C},h_0>0$, depending only on $T,\mathcal{O},M_1$, and $M_2$, such that for every $0<h\leq h_0$ and every $y_h^0\in L^2_{\mathcal F_0}(\Omega;H_h)$, the stochastic system \cref{sd:eq:stochastic-system} is exactly null controllable. That is, there exist controls $u_h\in L^2_{\mathbb F}(0,T;\mathcal U_h)$ and $v_h\in L^2_{\mathbb F}(0,T;H_h)$ such that the corresponding solution satisfies 
\begin{equation}
y_h(T)=0\quad\mathbb{P}\text{-a.s.}
\label{sd:eq:stochastic-exact-terminal}
\end{equation}
Moreover, these controls satisfy the mesh-independent cost estimate
\begin{equation}
\mathbb{E}\int_0^T
\bigl(\lVert u_h(t)\rVert_{\mathcal U_h}^2+\lVert v_h(t)\rVert_h^2\bigr)\,d t
\leq \mathcal{C}\mathbb{E}\lVert y_h^0\rVert_h^2.
\label{sd:eq:stochastic-exact-cost}
\end{equation}
\end{theorem}
\begin{proof}
By the Borel measurability of the minimum-norm control map $\mathcal K_h$ defined in \cref{sd:eq:minimal-control}, and the strong $\mathcal F_T$-measurability of $Q_h$ established in \Cref{sec:deterministic-exact}, the control $\widetilde u_h=\mathcal K_h(Q_h)y_h^0$ is strongly $\mathcal{F}_T$-measurable. The finite-dimensional solution map gives the same measurability for its trajectory $\widetilde y_h$ in $C([0,T];H_h)$. By \cref{sd:eq:Kh-bound}, the selected control is bounded uniformly in $h$. The energy estimate for its trajectory gives \cref{sd:eq:bridge-uniform-bound} with $R_h=0$. Applying \cref{sd:lem:mesh-lifting} yields \cref{sd:eq:stochastic-exact-terminal,sd:eq:stochastic-exact-cost}.
\end{proof}
\begin{remark}
For bounded sampled scalar coefficients, the relaxed  estimate
appears in \cite[Theorem~1.1]{Zhao2025}.
Lifting the controlled trajectories supplied by \cref{sd:prop:relaxed-control}
gives stochastic controls with an exponentially small terminal remainder.
In \cref{sd:thm:deterministic-exact}, a short-time deterministic correction
removes this remainder before lifting. Endpoint preservation then gives
exact stochastic null controllability.
\end{remark}
\begin{remark}\label{sd:rem:exact-scope}
The exterior-power comparison in \cref{sd:lem:determinant-comparison} exploits two key structural properties: the nearest-neighbour (tridiagonal) structure of the discrete Laplacian, which prevents the emergence of negative signs during basis expansions in the exterior algebra, and the diagonal form of $Q_h$, which preserves the diagonal structure in the compound matrices. Together, these properties guarantee the Metzler property of the compound generators. 
It does not require $Q_h$ to commute with $A_h$. The two adjacent observed
nodes used in \cref{sd:prop:short-observability} exist on every sufficiently
fine mesh because $\mathcal O$ is an interval. Extending the short-time
estimate to nonuniform meshes, higher dimensions, or non-diagonal drift
couplings requires a different argument. No diagonal structure is needed
for $C_h$: the unrestricted diffusion control realizes
$G_h-C_hy_h$ after projection.
\end{remark}

\section{Constrained controllability to terminal random targets}\label{sec:positive-parabolic}

In the previous two sections, the controls are unrestricted in sign and magnitude, and there is no state constraint. In this section, we address a more challenging scenario: reaching terminal states of uniformly bounded reference trajectories while enforcing a drift control in a fixed positive box and preserving nonnegativity of the state.

For this purpose, a long waiting time is employed. During the waiting time, the initial discrepancy decays exponentially due to the dissipativity of the parabolic semigroup. A positive margin in the reference control box then provides enough room to absorb a bounded null correction. On the deterministic side, this construction follows the constrained-control mechanisms in \cite{LoheacTrelatZuazua2017,PighinZuazua2018}. The internal-control extension is noted in \cite[Remark~1]{PighinZuazua2018}. Because the terminal reference data are random and only $\mathcal F_T$-measurable, the deterministic trajectory is generally not adapted. We will apply the general adapted lifting framework from \cref{sec:framework} to preserve the terminal endpoint, the pointwise control bounds, and the state nonnegativity.

\subsection{The constrained stochastic parabolic equation}

We retain the geometric setting and the spaces $H=L^2(G)$, $V=H_0^1(G)$, $\mathcal D=H^2(G)\cap H_0^1(G)$, and $U_1=L^2(G_0)$ introduced in \cref{sec:parabolic-setting}. Here $G$ is a bounded connected domain with $C^2$ boundary and $G_0\Subset G$ is nonempty and open. The operator $B:U_1\to H$ denotes the extension by zero.

In this section, we focus on the case where the principal operator is the standard Dirichlet Laplacian $A_D=-\Delta$ with $\operatorname{Dom}(A_D)=\mathcal D$. Let $S(t)=\mathrm{e}^{-tA_D}$ be the associated analytic $C_0$-semigroup. By the maximum principle, $S(t)$ is a positive semigroup. Furthermore, let $\lambda_1>0$ denote the first Dirichlet eigenvalue; then the semigroup satisfies the decay bound $\lVert S(t)\rVert_{\mathcal L(H)}\leq e^{-\lambda_1t}$ for $t\geq0$. We continue to use the trajectory space $\mathscr{X}(s,r)$ from \cref{eq:abstract-intersection-space} with the present spaces $H$ and $\mathcal D$.

We consider the control system
\begin{equation}\label{eq:positive-system}
\begin{cases}
d y=(-A_Dy+Bu)\,d t+(a_2y+v)\,d W(t) & 
t \in (0,T],\\[2pt]
y(0)=y_0,
\end{cases}
\end{equation}
where $a_2\in L^\infty_{\mathbb F}(0,T;L^\infty(G))$.
The drift control is localized in $G_0$, while the diffusion control acts
on all of $G$ and is unrestricted in sign.

The state constraint is the closed convex cone
$H_+:=\{h\in H:h\geq0\text{ a.e. in }G\}$.
For $0\leq \cM_1<\cM_2<\infty$, define the closed convex control box
\begin{equation*}
K_{\cM_1,\cM_2}:=\{q\in U_1:\cM_1\leq q(x)\leq \cM_2\text{ a.e. in }G_0\}.
\end{equation*}
The corresponding set of deterministic control paths on $(s,r)$ is
\begin{equation*}
\mathscr K_{\cM_1,\cM_2}(s,r)
:=\{q\in L^2(s,r;U_1):q(t)\in K_{\cM_1,\cM_2}\text{ for a.e. }t\in(s,r)\}.
\end{equation*}
The target will be the terminal value of a reference trajectory. Its data
may depend on all the information available at time $T$.

Fix $0<\nu<R<\infty$ and $\overline R_0>0$.
The reference initial state
$\overline y_0\in L^\infty_{\mathcal F_T}(\Omega;H)$ satisfies\vspace{-1mm}
\begin{equation}\label{eq:positive-ybar0}
\lVert\overline y_0(\omega)\rVert_H\leq\overline R_0
\quad\text{a.s.}
\end{equation}
The reference control satisfies\vspace{-1mm}
\begin{equation}\label{eq:rc}
\overline u\in
L^\infty\bigl(\Omega\times(0,T)\times G_0,
\mathcal F_T\otimes\mathcal B((0,T)\times G_0)\bigr),
\end{equation}
where $\mathcal B((0,T)\times G_0)$ is the Borel $\sigma$-field.
The reference control also satisfies\vspace{-1mm}
\begin{equation}
\nu\leq \overline u(\omega,t,x)\leq R
\quad\text{for a.e. }(\omega,t,x).
\label{eq:positive-ubar-box}\vspace{-1mm}
\end{equation}
The corresponding reference trajectory is\vspace{-1mm}
\begin{equation}\label{eq:positive-reference-trajectory}
\overline y(t)
=S(t)\overline y_0+
\int_0^tS(t-r)B\overline u(r)\,d r.\vspace{-1mm}
\end{equation}
The strict margin in \cref{eq:positive-ubar-box} absorbs a signed null
correction after a sufficiently long waiting time (the uniform bounds
on the reference data allow the waiting time to be deterministic).

\subsection{Samplewise deterministic bridge inside the control box}\label{sec:positive-bridge}

Following the samplewise deterministic reduction, we fix the sample point $\omega$ and construct a deterministic control that bridges the initial state and the terminal reference state, while satisfying the control box constraints.

The key mechanism is to exploit the strict positive margin $\nu>0$ in the reference control \eqref{eq:positive-ubar-box}. For a sufficiently long time horizon, the discrepancy between the initial state and the reference trajectory decays exponentially due to the dissipativity of the parabolic semigroup. This leaves a small residual at an intermediate time, which can be eliminated by a short-time null control that fits into the margin $\nu/2$. Consequently, the combined control remains within the enlarged box $[\nu/2, R+\nu/2]$.

To be precise, we first recall a standard deterministic null controllability result. For zero potential and a unit control interval, the results of \cite[Theorem~3.1 and Remark~3.1]{FernandezCaraZuazua2000} give, for every $h\in H$, a control $w\in L^\infty((0,1)\times G_0)$ such that
\begin{equation}\label{eq:positive-correction-system}
z'+A_Dz=Bw,\qquad z(0)=h,\quad z(1)=0,
\end{equation}
and
\begin{equation}\label{eq:positive-correction-cost}
\lVert w\rVert_{L^\infty((0,1)\times G_0)}\leq \mathcal{C}\lVert h\rVert_H.
\end{equation}
The constant $\mathcal C$ depends only on $G$ and $G_0$.

Fix $R_0,\overline R_0>0$ and $0<\nu<R$. Write $|G_0|$ for the Lebesgue
measure of $G_0$. With $\mathcal C$ fixed as in
\cref{eq:positive-correction-cost}, choose $T_+\geq2$, depending only on
$G,G_0,R_0,\overline R_0,\nu,R$, such that
\begin{equation}\label{eq:positive-Tplus}
\mathcal{C}\mathrm{e}^{-\lambda_1(T_+-2)}
\bigl(R_0+\overline R_0+R\lvert G_0\rvert^{1/2}\bigr)
\leq\frac\nu2.
\end{equation}
We use the reference control while the initial discrepancy decays, then
remove the remaining discrepancy during the final unit interval.
\begin{proposition}
\label{prop:positive-bridge}
Let $T\geq T_+$. Suppose $X\in V$ satisfies $\lVert X\rVert_H\leq R_0$, and let $(\overline y_0,\overline u)$ be deterministic reference data satisfying \cref{eq:positive-ybar0,eq:positive-ubar-box}. Then there exists a control $\widetilde u\in\mathscr K_{\nu/2,R+\nu/2}(1,T)$ such that the solution to 
\[
\widetilde y'+A_D\widetilde y=B\widetilde u, \qquad \widetilde y(1)=X
\]
belongs to $\mathscr X(1,T)\cap C([1,T];V)$ and reaches the terminal reference state, i.e., $\widetilde y(T)=\overline y(T)$.

Moreover, this control and the resulting trajectory satisfy the bounds
\begin{align}
&\lVert\widetilde u\rVert_{L^2(1,T;U_1)}^2
\leq (R+\nu/2)^2\lvert G_0\rvert(T-1), \label{eq:positive-bridge-control-bound} \\
&\lVert\widetilde y\rVert_{\mathscr X(1,T)}^2
\leq \mathcal C(T)\bigl(
\lVert X\rVert_V^2+(R+\nu/2)^2\lvert G_0\rvert(T-1)\bigr), \label{eq:positive-bridge-state-bound}
\end{align}
where $\mathcal C(T)$ depends only on $T$ and $G$. Finally, if $X\in H_+$, the trajectory preserves positivity, i.e., $\widetilde y(t)\in H_+$ for all $t\in[1,T]$.
\end{proposition}
\begin{proof}
The construction is carried out in two phases.

\textbf{Phase 1: Waiting and observation.}
On the interval $[1,T-1]$, we apply the reference control $\overline u$. The discrepancy at time $T-1$ is given by $h=S(T-2)(X-\overline y(1))$. By the bound on $\overline u$ and the exponential decay of the semigroup, we have
\begin{equation*}
\lVert h\rVert_H\leq e^{-\lambda_1(T-2)}
\bigl(R_0+\overline R_0+R\lvert G_0\rvert^{1/2}\bigr).
\end{equation*}

\textbf{Phase 2: Short-time null control.}
On the final unit interval $[T-1,T]$, we apply the short-time null control defined in \cref{eq:positive-correction-system,eq:positive-correction-cost} to drive the discrepancy $h$ to zero. By the choice of $T_+$ in \cref{eq:positive-Tplus}, this control satisfies $\lVert w\rVert_\infty\leq\mathcal C\lVert h\rVert_H\leq\nu/2$. 
Define $\widetilde u(t)=\overline u(t)+w(t-T+1)$ for $T-1<t<T$. Since this null control drives the discrepancy to zero, we have $\widetilde y(T)=\overline y(T)$. Furthermore, the combined control satisfies $\nu/2\leq\widetilde u\leq R+\nu/2$, which establishes \cref{eq:positive-bridge-control-bound}.

Finally, autonomous maximal $L^2$-regularity \cite[Theorem~4.1]{Dier2015}, together with the Dirichlet domain identifications on a $C^2$ domain, ensures that the trajectory belongs to the asserted space and satisfies \cref{eq:positive-bridge-state-bound} from $X\in V$ and $B\widetilde u\in L^2(1,T;H)$. The positivity of $S$ and $B$ yields the final assertion.
\end{proof}

\subsection{Constrained stochastic controllability to terminal random targets}\label{sec:positive-main}

We are now in a position to establish the main result of this section. In the previous subsection, we constructed a samplewise deterministic bridge that drives an initial state $X\in V$ with $\lVert X\rVert_H\leq R_0$ to the terminal reference target while satisfying the control constraints. However, because the terminal reference data are only $\mathcal{F}_T$-measurable, the selected deterministic control generally depends on future information and fails to be adapted.

In this subsection, we combine the measurable constrained selector from \cref{sec:measurable-common} and the adapted lifting from \cref{sec:projection}. The constrained selector makes the preceding deterministic construction measurable, and applying \cref{thm:abstract-lifting} preserves its terminal endpoint, control bounds, and state positivity.
\begin{theorem}
\label{thm:positive-main}
Let $T\geq T_+$, with $T_+$ chosen as in \cref{eq:positive-Tplus}, and let $y_0\in L^\infty_{\mathcal F_0}(\Omega;H)$ satisfy $\lVert y_0\rVert_{L^\infty(\Omega;H)}\leq R_0$. For reference data $(\overline y_0,\overline u)$ satisfying \cref{eq:positive-ybar0,eq:rc,eq:positive-ubar-box}, there exist controls $u\in L^2_{\mathbb F}(0,T;U_1)$ and $v\in L^2_{\mathbb F}(0,T;H)$ such that $u=0$ on the preparation interval $(0,1)$, and
\begin{equation*}
\nu/2\leq u(t,x)\leq R+\nu/2 \quad \text{for a.e. } (\omega,t,x)\in\Omega\times(1,T)\times G_0.
\end{equation*} 
Furthermore, the solution of \cref{eq:positive-system} reaches the terminal random target
$
y(T)=\overline y(T)$, $\mathbb P$-a.s. 
The state and controls satisfy the energy estimate
\begin{equation}\label{eq:positive-cost}
\mathbb{E}\sup_{0\leq t\leq T}\lVert y(t)\rVert_H^2
+\mathbb{E}\int_0^T\bigl(\lVert u(t)\rVert_{U_1}^2+\lVert v(t)\rVert_H^2\bigr)\,d t
\leq \mathcal{C}
\bigl(R_0^2+\overline R_0^2+R^2\lvert G_0\rvert+1\bigr),
\end{equation}
where $\mathcal{C}$ depends only on $T,G,G_0,\nu,R$, and $\lVert a_2\rVert_{L^\infty(\Omega\times(0,T)\times G)}$. Finally, if $y_0\in H_+$ almost surely, the corresponding state preserves positivity, i.e.,
\begin{equation*}
\mathbb{P}\bigl(y(t)\in H_+\text{ for every }t\in[0,T]\bigr)=1.
\end{equation*}
\end{theorem}

\begin{proof}
\textbf{Step 1. Prepare the switching state.}
On $[0,1]$, take $u=0$ and $v=-a_2S(t)y_0$.
This cancels the noise, so the preparation trajectory is $X(t)=S(t)y_0$. Analytic smoothing gives
$X(1)\in L^2_{\mathcal F_1}(\Omega;V)$ and
\begin{equation}
\mathbb{E}\lVert X(1)\rVert_V^2
\leq \mathcal{C}\mathbb{E}\lVert y_0\rVert_H^2.
\label{eq:positive-prepared-state}
\end{equation}
By the contraction property of the semigroup, $\lVert X(1)\rVert_H\leq R_0$, and the positivity of the semigroup preserves $H_+$ on $[0,1]$ for nonnegative initial data.

\smallskip

\textbf{Step 2. Select a constrained trajectory in $\mathscr{X}(1,T)$.}
We fix the interval $(1,T)$ and define the bounded endpoint operator $\Lambda_T:L^2(1,T;U_1)\to H$ by $\Lambda_Tq=\int_1^T S(T-r)Bq(r)\,d r$. 
By \cref{eq:positive-reference-trajectory}, the residual $g=\overline y(T)-S(T-1)X(1)$ is strongly $\mathcal F_T$-measurable. The deterministic bridge from \cref{prop:positive-bridge} guarantees that $g(\omega)\in \Lambda_T(\mathscr K_T)$ almost surely, where $\mathscr K_T=\mathscr K_{\nu/2,R+\nu/2}(1,T)$.
By \cref{prop:constrained-selector}, there is a strongly $\mathcal F_T$-measurable map $\widetilde u:\Omega\to L^2(1,T;U_1)$ such that $\widetilde u(\omega)\in\mathscr K_T$ and $\Lambda_T\widetilde u(\omega)=g(\omega)$ almost surely.
Since $\mathscr K_T$ is bounded in $L^2(1,T;U_1)$, $\widetilde u$ is square integrable. The deterministic solution map is continuous by  maximal $L^2$-regularity \cite[Theorem~4.1]{Dier2015}. Hence it gives a measurable trajectory $\widetilde y\in L^2_{\mathcal F_T}(\Omega;\mathscr{X}(1,T))$ satisfying \cref{eq:positive-bridge-state-bound} and the terminal condition $\widetilde y(T)=\overline y(T)$. Positivity of $S$ and $B$ ensures that $\widetilde y(t)\in H_+$ for all $t\in[1,T]$ if $X(1)\in H_+$.

\textbf{Step 3. Adapt the trajectory and concatenate.}
We now apply the adapted lifting \cref{thm:abstract-lifting} to the pair $(\widetilde y, \widetilde u)$ on $(1,T)$.
Set $\mathcal A_0=-A_D$, $N=0$,
$C(t,\omega)z=a_2(t,\omega,\cdot)z$, and $f=0$.
We enforce the fixed control constraint $K_{U_1}=K_{\nu/2,R+\nu/2}\subset U_1$.
The lifting produces an adapted state $y_+$ and controls
$u_+=\Pi_{U_1}\widetilde u$ and $v_+=\Sigma-a_2y_+$,
preserving the endpoint $\overline y(T)$ and the control box.
If $y_0\in H_+$ almost surely, then the preparation and the
deterministic bridge preserve nonnegativity. Applying the lifting
with $K_H=H_+$ preserves this state constraint as well.

Finally, we concatenate the preparation trajectory on $[0,1]$ with the lifted trajectory on $[1,T]$. Since the states agree at the deterministic switching time $t=1$, the concatenated process solves \cref{eq:positive-system} on $[0,T]$. For this Laplacian model, standard variational well-posedness \cite[Proposition~9.3]{Lue2021a} gives 
\begin{align}\label{eq:stochastic-energy}
\mathbb{E}\sup_{0\leq t\leq T}\lVert y(t)\rVert_H^2
+\mathbb{E}\int_0^T\lVert y(t)\rVert_V^2\,d t
\leq \mathcal{C}\Bigl(\mathbb{E}\lVert y_0\rVert_H^2
+\mathbb{E}\int_0^T
\bigl(\lVert u(t)\rVert_{U_1}^2+\lVert v(t)\rVert_H^2\bigr)\,d t\Bigr),
\end{align}
with a constant depending only on $T,G,\lVert a_2\rVert_{L^\infty(\Omega\times(0,T)\times G)}$, and gives uniqueness by applying the same estimate to differences.
The bounds \cref{eq:positive-prepared-state,eq:positive-bridge-state-bound} and the lifting estimates \cref{eq:abstract-U-estimate,eq:abstract-V-estimate} control $u$ and $v$. Applying \cref{eq:stochastic-energy} then gives the total cost estimate \cref{eq:positive-cost}. When $y_0\in H_+$ almost surely, the state constraint holds at
all times almost surely, by continuity and closedness of $H_+$,
as in \cref{thm:abstract-lifting}.
\end{proof}

\section{Exact controllability of stochastic wave equations with internal controls}\label{sec:application-wave}

We next apply the lifting method to a controlled stochastic wave equation. Deterministic
internal observability gives controlled trajectories in the energy space.
A change of variables provides the time regularity needed to apply
\cref{thm:abstract-lifting} with a zero generator.

\subsection{Formulation of the problem and the controllability result}

Recall that $G\subset\mathbb R^n$ is a bounded domain with a
$C^2$ boundary. Let $G_0\subset G$ be open. Fix
$x_0\in\mathbb R^n\setminus\overline G$ and $\delta>0$. Writing $\nu$
for the outward unit normal, set
\begin{equation*}
\Gamma_0=\{x\in\partial G:(x-x_0)\cdot\nu(x)>0\}.
\end{equation*}
We assume that the control region contains an interior neighborhood of
this part of the boundary:
\begin{equation}\label{eq:wave-control-geometry}
G_\delta:=\{x\in G:\operatorname{dist}(x,\Gamma_0)<\delta\}\subset G_0.
\end{equation}
Decreasing $\delta$ if necessary, we assume that $G\setminus G_\delta$
is nonempty. Define
\begin{equation}\label{eq:wave-control-time}
T_0=2\Bigl(\sup_{x\in G\setminus G_\delta}|x-x_0|^2
-\inf_{x\in G}|x-x_0|^2\Bigr)^{1/2}.
\end{equation}
By \cite[Theorem~1.4 and Remark~1.7]{FuLiaoLu2023}, deterministic
internal observability holds for every $T>T_0$ and every bounded
time-dependent potential. The factor $2$ in \cref{eq:wave-control-time}
accounts for the symmetric time interval used in that reference.

We take $\mathfrak K=\mathbb R$ and set
$H=H_0^1(G)\times L^2(G)$ and $U_1=L^2(G_0)$. The energy norm is
\begin{equation*}
\lVert(q,p)\rVert_H^2=\lVert\nabla q\rVert_{L^2(G;\mathbb R^n)}^2
+\lVert p\rVert_{L^2(G)}^2.
\end{equation*}
Let $a,b_2,b_3\in L^\infty_{\mathbb F}(0,T;L^\infty(G))$ and $b_1\in L^\infty_{\mathbb F}(0,T;W^{1,\infty}(G))$.
Set
\begin{equation}\label{eq:wave-coefficient-bound}
\cM_3:= |a|_{L^\infty_{\mathbb F}(0,T;L^\infty(G))} + |b_1|_{L^\infty_{\mathbb F}(0,T;W^{1,\infty}(G))} + |b_2|_{L^\infty_{\mathbb F}(0,T;L^\infty(G))} + |b_3|_{L^\infty_{\mathbb F}(0,T;L^\infty(G))}.
\end{equation}
With $Q=(0,T)\times G$ and $\Sigma=(0,T)\times\partial G$, consider
\begin{equation}\label{eq:wave-system}
\begin{cases}
d y=\widehat y\,d t+(b_1y+v_1)\,d W(t) &\text{in }Q,\\
d\widehat y=(\Delta y+ay+\mathbf 1_{G_0}u)\,d t
+(b_2y+b_3\widehat y+v_2)\,d W(t) &\text{in }Q,\\
y=0 &\text{on }\Sigma,\\
(y(0),\widehat y(0))=(y_0,\widehat y_0) &\text{in }G.
\end{cases}
\end{equation}
Here $y$ is the displacement and $\widehat y$ is the velocity.
The controls satisfy $u\in L^2_{\mathbb F}(0,T;U_1)$ and
$v=(v_1,v_2)\in L^2_{\mathbb F}(0,T;H)$.  

The Dirichlet wave operator is
\begin{equation*}
\mathcal A_0(q,p)=(p,\Delta q),\qquad
\operatorname{Dom}(\mathcal A_0)
=(H^2(G)\cap H_0^1(G))\times H_0^1(G).
\end{equation*}
It generates a unitary group $\{S_0(t)\}_{t\in\mathbb R}$ on $H$.
By \cite[Theorem 3.20]{Lue2021a},
\cref{eq:wave-system} has a unique mild solution
$X=(y,\widehat y)\in L^2_{\mathbb F}(\Omega;C([0,T];H))$ for every
$X_0=(y_0,\widehat y_0)\in L^2_{\mathcal F_0}(\Omega;H)$ and every
control pair $(u,v)$ in the stated spaces.
\begin{theorem}\label{thm:wave-exact}
Assume the geometric condition \cref{eq:wave-control-geometry} and the coefficient bound \cref{eq:wave-coefficient-bound}, and let $T_0$ be given by \cref{eq:wave-control-time}. For every $T>T_0$, the control system \cref{eq:wave-system} is exactly controllable.  
More precisely, for any initial state $X_0\in L^2_{\mathcal F_0}(\Omega;H)$ and any target $\xi\in L^2_{\mathcal F_T}(\Omega;H)$, there exist controls $u\in L^2_{\mathbb F}(0,T;U_1)$ and $v\in L^2_{\mathbb F}(0,T;H)$ such that the corresponding solution satisfies
$X(T)=\xi$, $\mathbb P$-a.s. 
Furthermore, these controls satisfy the energy and cost estimate
\begin{equation}\label{eq:wave-control-cost}
\mathbb E\sup_{0\leq t\leq T}\lVert X(t)\rVert_H^2
+\mathbb E\int_0^T
\bigl(\lVert u(t)\rVert_{U_1}^2+\lVert v(t)\rVert_H^2\bigr)\,d t
\leq\mathcal C
\bigl(\mathbb E\lVert X_0\rVert_H^2+\mathbb E\lVert\xi\rVert_H^2\bigr),
\end{equation}
where $\mathcal C$ depends only on the fixed geometry, $T$ and $\cM_3$.
\end{theorem}

\subsection{Proof of \texorpdfstring{\cref{thm:wave-exact}}{the wave controllability theorem}}

\begin{proof}
\textbf{Step 1. Select a measurable deterministic controlled trajectory.}
Fix a sample point outside a null set on which the coefficient bound may
fail. Apply \cite[Theorem~1.4 and Remark~1.7]{FuLiaoLu2023} in the
scalar case with the identity principal matrix and $p=\infty$, translating
the time interval $(-T/2,T/2)$ to $(0,T)$. Its observation region is
contained in $(0,T)\times G_\delta\subset(0,T)\times G_0$.
Thus, for $T>T_0$ with $T_0$ given by \cref{eq:wave-control-time},
every weak solution of
$z_{tt}-\Delta z=a(\omega)z$ with homogeneous Dirichlet boundary data
and initial data in $L^2(G)\times H^{-1}(G)$
satisfies
\begin{equation*}
\lVert z(0)\rVert_{L^2(G)}^2
+\lVert z_t(0)\rVert_{H^{-1}(G)}^2
\leq\mathcal C\int_0^T\int_{G_0}|z(t,x)|^2\,d x\,d t.
\end{equation*}
The constant is uniform over potentials bounded by $\cM_3$. Wave duality
and the energy estimate therefore yield a deterministic controlled pair
with $\widetilde X=(\widetilde y,\partial_t\widetilde y)$ satisfying
\begin{equation}\label{eq:wave-deterministic-system}
\begin{cases}
\partial_{tt}\widetilde y
=\Delta\widetilde y+a(\omega)\widetilde y
+\mathbf 1_{G_0}\widetilde u &\text{in }Q,\\
\widetilde y=0 &\text{on }\Sigma,\\
\widetilde X(0)=X_0(\omega),\qquad
\widetilde X(T)=\xi(\omega).
\end{cases}
\end{equation}
The minimum-norm control and its trajectory satisfy
\begin{equation}\label{eq:wave-deterministic-bound}
\sup_{0\leq t\leq T}\lVert\widetilde X(t)\rVert_H^2
+\int_0^T\lVert\widetilde u(t)\rVert_{U_1}^2\,d t
\leq\mathcal C
\bigl(\lVert X_0(\omega)\rVert_H^2+\lVert\xi(\omega)\rVert_H^2\bigr).
\end{equation}
To ensure that this samplewise selection is measurable, let $\Lambda_T(\omega)$ be the endpoint
map with zero initial data and let $X_T^0(\omega)$ be the uncontrolled
terminal state. By applying Picard iteration to the mild formulation of the wave equation, the solution is obtained as the limit of successive approximations. Since the potential is predictable and the control operator is deterministic and the semigroup is deterministic, each iteration is strongly $\mathcal F_T$-measurable. The uniform bound on the potential and the unitarity of the wave group give the usual factorial bound on successive Picard differences, so the sequence converges in $L^2(\Omega; C([0,T];H))$. Consequently, the terminal endpoint map $\Lambda_T$ is strong-operator $\mathcal F_T$-measurable and uniformly bounded, and $X_T^0\in L^2_{\mathcal F_T}(\Omega;H)$. 
Apply \cref{thm:random-Tikhonov} to
$L=\Lambda_T$ and $g=\xi-X_T^0$.
The bound in \cref{eq:wave-deterministic-bound} gives a measurable
minimum-norm control
$\widetilde u\in L^2_{\mathcal F_T}(\Omega;L^2(0,T;U_1))$.
The wave equation \cref{eq:wave-deterministic-system} then gives
$\widetilde X\in L^2_{\mathcal F_T}(\Omega;C([0,T];H))$.

\smallskip
\textbf{Step 2. Apply the lifting theorem after a change of variables.}
Write $N(t)(q,p)=(0,a(t)q)$ and
$Bu=(0,\mathbf 1_{G_0}u)$. The diffusion operator is
$C(t)(q,p)=(b_1(t)q,b_2(t)q+b_3(t)p)$.
By \cref{eq:wave-coefficient-bound}, $N$ and $C$ are strongly
predictable and uniformly bounded on $H$, and
$B\in\mathcal L(U_1;H)$.

Set $\widetilde Y(t)=S_0(-t)\widetilde X(t)$.
The mild form of \cref{eq:wave-deterministic-system} gives
\begin{equation*}
\widetilde Y(t)=X_0+\int_0^t
S_0(-r)\bigl[N(r)\widetilde X(r)+B\widetilde u(r)\bigr]\,d r.
\end{equation*}
Thus $\widetilde Y\in L^2_{\mathcal F_T}(\Omega;H^1(0,T;H))$,
and its derivative satisfies
\begin{equation}\label{eq:wave-transformed-bound}
\mathbb E\int_0^T\lVert\dot{\widetilde Y}(t)\rVert_H^2\,d t
\leq\mathcal C\mathbb E\int_0^T
\bigl(\lVert\widetilde X(t)\rVert_H^2
+\lVert\widetilde u(t)\rVert_{U_1}^2\bigr)\,d t.
\end{equation}
Define $\widehat N(t)=S_0(-t)N(t)S_0(t)$,
$\widehat B(t)=S_0(-t)B$  and
$\widehat C(t)=S_0(-t)C(t)S_0(t)$.
These operators are strongly predictable and uniformly bounded.
We apply \cref{thm:abstract-lifting} with generator zero, common domain
$H$, and coefficients $\widehat N,\widehat B,\widehat C$.
It gives $u=\Pi_{U_1}\widetilde u$, a diffusion control
$\widehat v\in L^2_{\mathbb F}(0,T;H)$, and the continuous adapted
process $Y(t)=\mathbb E[\widetilde Y(t)\mid\mathcal F_t]$ satisfying
\begin{equation}\label{eq:wave-transformed-stochastic}
\begin{cases}
dY=(\widehat N Y+\widehat B u)\,d t
+(\widehat C Y+\widehat v)\,d W(t),\\
Y(0)=X_0,\qquad Y(T)=S_0(-T)\xi.
\end{cases}
\end{equation}

\smallskip
\textbf{Step 3. Recover the wave equation and estimate the controls.}
Set $X(t)=S_0(t)Y(t)$ and $v(t)=S_0(t)\widehat v(t)$.
Multiplying the integral form of
\cref{eq:wave-transformed-stochastic} by $S_0(t)$ yields
\begin{equation*}
X(t)= S_0(t)X_0
+\int_0^tS_0(t-r)\bigl[N(r)X(r)+Bu(r)\bigr]\,d r +\int_0^tS_0(t-r)\bigl[C(r)X(r)+v(r)\bigr]\,d W(r).
\end{equation*}
This is the mild form of \cref{eq:wave-system}, and the endpoint
identities in \cref{eq:wave-transformed-stochastic} give
$X(T)=\xi$, $\mathbb P$-a.s. 

Since $S_0(t)$ is unitary, \cref{eq:abstract-U-estimate,eq:abstract-V-estimate,eq:wave-transformed-bound,eq:wave-deterministic-bound} imply
\begin{equation*}
\mathbb E\int_0^T
\bigl(\lVert u(t)\rVert_{U_1}^2+\lVert v(t)\rVert_H^2\bigr)\,d t
\leq\mathcal C
\bigl(\mathbb E\lVert X_0\rVert_H^2+\mathbb E\lVert\xi\rVert_H^2\bigr).
\end{equation*}
Moreover, $X(t)=\mathbb E[\widetilde X(t)\mid\mathcal F_t]$ for each deterministic $t$. Conditional Jensen's inequality, Doob's maximal inequality, and the pathwise continuity of $X$ give
\begin{equation*}
\mathbb E\sup_{0\leq t\leq T}\lVert X(t)\rVert_H^2
\leq4\mathbb E\sup_{0\leq t\leq T}\lVert\widetilde X(t)\rVert_H^2.
\end{equation*}
Combining this with \cref{eq:wave-deterministic-bound} proves
\cref{eq:wave-control-cost}.
\end{proof}

\begin{remark}\label{rem:wave-diffusion-channels}
The pair $(v_1,v_2)$ realizes the martingale correction in the full energy
space $H$. A diffusion control acting only in the momentum equation has
range $\{0\}\times L^2(G)$ and does not satisfy the right-inverse
assumption of \cref{cor:endpoint-transfer}.
\end{remark}

\appendix

\section{Proofs for semi-discrete spectral and determinant estimates}\label{app:semidiscrete-technical}

In this appendix, we provide the detailed proofs for the endpoint observability estimate used for each grid block in \cref{sd:prop:short-observability}. The proof is divided into three main steps:
\begin{enumerate}[label=\textup{(\roman*)},font=\upshape]
\item In \cref{sd:lem:finite-biorthogonal}, we construct a finite biorthogonal family in time to obtain the zero-potential endpoint observability estimate.
\item In \cref{sd:lem:determinant-comparison}, we exploit the positivity of the evolution operators on exterior power spaces, specifically, the Metzler property of the compound matrices, to compare the mixed observation determinants for bounded diagonal potentials with the zero-potential case.
\item In \cref{sd:prop:robust-endpoint}, we apply Andr\'eief's identity to compare the coefficients of the Gram determinant, yielding the explicit dependence on the block size, potential bound, and control time required in \cref{sd:eq:robust-constant}.
\end{enumerate}
Throughout the proofs of \cref{sd:lem:finite-biorthogonal} and \cref{sd:lem:unperturbed-endpoint}, fix $0<L_-<L_+$ and let $n\geq1$, $h>0$, and $\delta$ satisfy $L_-\leq(n+1)h\leq L_+$ and $0<\delta\leq1$.
The eigenvalues of $A_{h,n}$ are
\begin{equation*}
\lambda_{k,h}=\frac4{h^2}\sin^2\biggl(\frac{k\pi}{2(n+1)}\biggr),
\qquad 1\leq k\leq n,
\end{equation*}
where the block size $n$ is suppressed in the notation.

The following lemma provides the standard construction of a finite biorthogonal family for the exponents $-\lambda_{k,h}$, with an explicit dependence on the short time $\delta$.
\begin{lemma}
\label{sd:lem:finite-biorthogonal}
There is $\mathcal C>0$, depending only on $L_-,L_+$, such that for every
admissible $n,h,\delta$ there are real functions
$\{\theta_{k,h}^{\delta}\}_{k=1}^n\subset L^2(0,\delta)$ satisfying
\begin{equation}
\int_0^\delta
\theta_{k,h}^{\delta}(t)e^{-\lambda_{\ell,h}t}\,d t
=\delta_{k\ell},\qquad 1\leq k,\ell\leq n,
\label{sd:eq:finite-biorthogonality}
\end{equation}
where $\delta_{k\ell}$ is the Kronecker symbol. Their norms satisfy
\begin{equation}
\lVert \theta_{k,h}^{\delta}\rVert_{L^2(0,\delta)}^2
\leq \mathcal{C}\exp\biggl(
\frac{\mathcal{C}}{\delta}+\frac{\delta}{2}\lambda_{k,h}
\biggr).
\label{sd:eq:finite-biorthogonal-bound}
\end{equation}
\end{lemma}
\begin{proof}
\textbf{Step 1. Interpolate the finite spectrum.}
Put $\vartheta_k=k\pi/(n+1)$. The monic characteristic polynomial is
$p_{h,n}(\lambda)=(-1)^nh^{-2n}U_n(1-h^2\lambda/2)$, where
$U_n$ is the Chebyshev polynomial of the second kind, given by
$U_n(\cos\vartheta)=\sin((n+1)\vartheta)/\sin\vartheta$.
Thus $|U_n'(\cos\vartheta_k)|=(n+1)/\sin^2\vartheta_k$ and
$\prod\limits_{j=1}^n\lambda_{j,h}=(n+1)h^{-2n}$, which give
\begin{equation*}
D_{k,h}:=\prod_{j\ne k}
\biggl|1-\frac{\lambda_{k,h}}{\lambda_{j,h}}\biggr|
=\frac{\lambda_{k,h}|p_{h,n}'(\lambda_{k,h})|}
{\prod_{j=1}^n\lambda_{j,h}}
=\frac1{2\cos^2(\vartheta_k/2)}\geq\frac12.
\end{equation*}
The interpolation polynomial
\begin{equation*}
P_{k,h}(z)=\prod_{j\ne k}
\frac{1+iz/\lambda_{j,h}}{1-\lambda_{k,h}/\lambda_{j,h}}
\end{equation*}
satisfies $P_{k,h}(i\lambda_{\ell,h})=\delta_{k\ell}$.
Using $\lambda_{j,h}\geq4j^2/L_+^2$ and comparing the sum with an integral, we have
\begin{equation*}
\sum_{j\geq1}\log\biggl(1+\frac{R^4}{j^4}\biggr)
\leq\int_0^\infty\log\biggl(1+\frac{R^4}{t^4}\biggr)\,d t
=\mathcal{C}R,\qquad R>0,
\end{equation*}
which yields
\begin{equation}\label{sd:eq:finite-interpolant-growth}
|P_{k,h}(x)|\leq2\exp(\mathcal{C}\sqrt{|x|}),\qquad x\in\mathbb{R},
\end{equation}
where $\mathcal{C}$ depends only on $L_+$.

\textbf{Step 2. Apply the multiplier estimate.}
Choose a sufficiently small universal constant $\varepsilon>0$.
For $\beta,\nu>0$ with $\beta\nu=(\pi+\varepsilon)^2$, define
\begin{equation*}
H_{\beta,\nu}(z)
=\mathcal{C}_\nu\int_{-1}^1
\exp\biggl(-\frac{\nu}{1-t^2}\biggr)e^{-i\beta tz}\,d t,
\end{equation*}
where $\mathcal C_\nu$ normalizes the nonnegative density to have integral one:
\begin{equation*}
\mathcal C_\nu^{-1}=\int_{-1}^1\exp\biggl(-\frac{\nu}{1-t^2}\biggr)\,d t.
\end{equation*}
The multiplier $H_{\beta,\nu}$ has exponential type at most $\beta$ and, by
\cite[Lemma~2.3]{Lissy2014} with $\alpha=2$, satisfies
\begin{equation}\label{sd:eq:multiplier-real-bound}
|H_{\beta,\nu}(x)|
\leq \mathcal{C}\sqrt{\nu+1}
\exp\biggl(\frac{3\nu}{4}-(\pi+\varepsilon/2)\sqrt{|x|}\biggr).
\end{equation}
With $\mathcal{C}$ fixed as in \cref{sd:eq:finite-interpolant-growth}, choose
$b\geq1$ such that $(\pi+\varepsilon/2)\sqrt b>\mathcal{C}+1$.
Set $\beta=\delta/(4b)$, $\nu=(\pi+\varepsilon)^2/\beta$,
$M_\delta(z)=H_{\beta,\nu}(bz)$, and $\tau=\delta/4$.
Then $M_\delta$ has type at most $\tau$ and $\nu\leq \mathcal{C}/\delta$.
Combining \cref{sd:eq:finite-interpolant-growth,sd:eq:multiplier-real-bound} gives
\begin{equation}\label{sd:eq:interpolant-product-bound}
|P_{k,h}(x)M_\delta(x)|\leq \mathcal{C}e^{\mathcal{C}/\delta}e^{-\sqrt{|x|}}.
\end{equation}
Since the defining density is nonnegative, even, and has total mass one,
\begin{equation}\label{sd:eq:multiplier-imaginary-bound}
M_\delta(i\lambda)
=\mathcal{C}_\nu\int_{-1}^1e^{-\nu/(1-t^2)}
\cosh(\tau\lambda t)\,d t\geq1,\qquad\lambda\geq0.
\end{equation}

\textbf{Step 3. Translate the Paley--Wiener interpolant.}
We suppress the dependence on $\delta$ in $F_{k,h}$ and $g_{k,h}$. Set
\begin{equation*}
F_{k,h}(z)=P_{k,h}(z)\frac{M_\delta(z)}{M_\delta(i\lambda_{k,h})}.
\end{equation*}
This entire function has type at most $\tau$, satisfies
$F_{k,h}(i\lambda_{\ell,h})=\delta_{k\ell}$ and
$F_{k,h}(-x)=\overline{F_{k,h}(x)}$, and has squared $L^2(\mathbb R)$ norm
at most $\mathcal C e^{\mathcal C/\delta}$ by
\cref{sd:eq:interpolant-product-bound,sd:eq:multiplier-imaginary-bound}.
The $L^2$ Paley--Wiener theorem therefore supplies a real $g_{k,h}$ supported
in $[-\tau,\tau]$ with $\lVert g_{k,h}\rVert_{L^2(\mathbb R)}^2
\leq\mathcal C e^{\mathcal C/\delta}$ and
\begin{equation*}
F_{k,h}(z)=\int_{-\tau}^{\tau}g_{k,h}(r)e^{-izr}\,d r.
\end{equation*}
Set $\theta_{k,h}^{\delta}(t)
=e^{\tau\lambda_{k,h}}g_{k,h}(\tau-t)$ for $0<t<2\tau$, and zero on
$(2\tau,\delta)$. Then
\begin{equation*}
\int_0^\delta\theta_{k,h}^{\delta}(t)e^{-\lambda_{\ell,h}t}\,d t
=e^{(\lambda_{k,h}-\lambda_{\ell,h})\tau}
F_{k,h}(i\lambda_{\ell,h})=\delta_{k\ell}.
\end{equation*}
The translation gives the required norm bound:
\begin{equation*}
\lVert \theta_{k,h}^{\delta}\rVert_{L^2(0,\delta)}^2
=e^{2\tau\lambda_{k,h}}\lVert g_{k,h}\rVert_{L^2(\mathbb R)}^2
\leq \mathcal{C}\exp\biggl(\frac{\mathcal{C}}{\delta}+\frac\delta2\lambda_{k,h}\biggr).
\end{equation*}
The choices of $\varepsilon$ and $b$ are independent of $h,n,k,\delta$,
which gives the stated uniformity.
\end{proof}

The endpoint amplitudes convert this spectral bound into observability.

\begin{corollary}
\label{sd:lem:unperturbed-endpoint}
Under the preceding length and time restrictions, there is a constant
$\mathcal C>0$, depending only on $L_-,L_+$, such that every solution of
$z'+A_{h,n}z=0$ satisfies
\begin{equation}\label{sd:eq:unperturbed-endpoint}
\bigl\lVert z(\delta)\bigr\rVert_{h,n}^2
\leq \mathcal{C}h^{-4}\exp\biggl(\frac{\mathcal{C}}{\delta}\biggr)
\int_0^\delta h\bigl\lvert z_n(t)\bigr\rvert^2\,d t.
\end{equation}
\end{corollary}

\begin{proof}
Set $L=(n+1)h$. The weighted orthonormal eigenvectors are
$\phi_{k,h}(j)=\sqrt{2/L}\sin(jk\pi/(n+1))$. At the observed endpoint,
$b_{k,h}=\sqrt h\,\phi_{k,h}(n)$ satisfies
\begin{equation}\label{sd:eq:endpoint-mode-lower-bound}
|b_{k,h}|=\sqrt{\frac{2h}{L}}\sin\biggl(\frac{k\pi}{n+1}\biggr)
\geq\frac{2\sqrt2}{L_+^{3/2}}h^{3/2}.
\end{equation}
Write $z(0)=\sum\limits_{k=1}^na_k\phi_{k,h}$ and
$y(t)=\sqrt h\,z_n(t)=\sum\limits_{k=1}^na_kb_{k,h}e^{-\lambda_{k,h}t}$.
By \cref{sd:eq:finite-biorthogonality},
$a_kb_{k,h}=\int_0^\delta y(t)\theta_{k,h}^{\delta}(t)\,d t$,
so the Cauchy--Schwarz inequality yields
\begin{equation*}
\lVert z(\delta)\rVert_{h,n}^2
\leq\lVert y\rVert_{L^2(0,\delta)}^2
\sum_{k=1}^n|b_{k,h}|^{-2}e^{-2\delta\lambda_{k,h}}
\lVert \theta_{k,h}^{\delta}\rVert_{L^2(0,\delta)}^2
\leq\mathcal C h^{-4}e^{\mathcal C/\delta}\lVert y\rVert_{L^2(0,\delta)}^2.
\end{equation*}
Here we used \cref{sd:eq:finite-biorthogonal-bound,sd:eq:endpoint-mode-lower-bound}
and $n\leq L_+/h$.
This proves \cref{sd:eq:unperturbed-endpoint}.
\end{proof}

For the remaining estimates, fix $M\geq0$, $h>0$, $n\geq1$, and
$0<\delta\leq1$. We use Euclidean coordinates and write $*$ for transpose.
The common weight $h$ in the state norm and endpoint
observation cancels in the observability quotient.
Let $q_j\in L^\infty(0,\delta;\mathbb R)$ satisfy
$\lVert q_j\rVert_{L^\infty(0,\delta)}\leq M$, and set
$Q_n(t)=\operatorname{diag}(q_1(t),\ldots,q_n(t))$.
Let $S_Q(t,s)$ denote the evolution operator for the system
\begin{equation*}
z'+A_{h,n}z=Q_n(t)z.
\end{equation*}
Write $S_0(t)=e^{-tA_{h,n}}$ for the semigroup with zero potential.
Define the observation Gramian by
\begin{equation}
W_Q=\int_0^\delta S_Q(t,0)^*e_ne_n^*S_Q(t,0)\,d t.
\label{sd:eq:tail-gramians}
\end{equation}
The terminal energy matrix is
\begin{equation}\label{sd:eq:tail-terminal-energy}
H_Q=S_Q(\delta,0)^*S_Q(\delta,0).
\end{equation}

To compare $W_Q$ and $H_Q$, we use determinants whose rows are endpoint
observations at distinct times or state observations at the terminal time.
The nearest-neighbour structure makes their absolute values comparable
despite the time dependence of $Q_n$.

\begin{lemma}
\label{sd:lem:determinant-comparison}
Let
$
\Gamma_\delta=
\bigl(\{n\}\times(0,\delta)\bigr)
\cup
\bigl(\{1,\ldots,n\}\times\{\delta\}\bigr),
$
and set $K_Q((i,t),j)=e_i^*S_Q(t,0)e_j$. For an $n$-tuple
$\boldsymbol\sigma=(\sigma_1,\ldots,\sigma_n)\in\Gamma_\delta^n$, write
$
D_Q(\boldsymbol\sigma)
:=\bigl\lvert\det[K_Q(\sigma_a,j)]_{a,j=1}^n\bigr\rvert.
$
Then
\begin{equation}\label{sd:eq:determinant-comparison}
e^{-nM\delta}D_0(\boldsymbol\sigma)
\leq D_Q(\boldsymbol\sigma)
\leq e^{nM\delta}D_0(\boldsymbol\sigma).
\end{equation}
\end{lemma}

\begin{proof}
\textbf{Step 1. Compare the exterior evolutions.}
Let $V=\mathbb R^n$, set $E_n=e_1\wedge\cdots\wedge e_n$, and equip
each exterior power $\bigwedge\limits^kV$ with the ordered orthonormal basis
\begin{equation*}
e_I=e_{i_1}\wedge\cdots\wedge e_{i_k},
\qquad I=(i_1<\cdots<i_k).
\end{equation*}
Put $L_0=-A_{h,n}$ and $L_Q(t)=L_0+Q_n(t)$.  The evolution induced by
$S_Q(t,s)$ on $\bigwedge\limits^kV$ is $\bigwedge\limits^kS_Q(t,s)$, and its generator
is the additive compound
\begin{equation*}
L_Q^{[k]}(t)=L_0^{[k]}+Q_n^{[k]}(t).
\end{equation*}
Here $L^{[k]}$ acts on a decomposable wedge by applying $L$ to one
factor at a time and summing over the $k$ factors.
Since $Q_n$ is diagonal and $\lvert q_j\rvert\leq M$,
\begin{equation*}
Q_n^{[k]}(t)e_I
=\Bigl(\sum_{i\in I}q_i(t)\Bigr)e_I.
\end{equation*}
Thus each diagonal entry of $Q_n^{[k]}$ lies in $[-kM,kM]$.

The matrix $L_0$ is tridiagonal with positive off-diagonal entries.  Every
off-diagonal transition in $L_0^{[k]}$ replaces one index in $I$ by a
nearest neighbour.  Whenever the resulting wedge is nonzero, this
nearest-neighbour replacement does not cross another index, so no negative
reordering sign is created. Hence $L_0^{[k]}$ is Metzler in the ordered
wedge basis, that is, it has nonnegative off-diagonal entries.
The comparison theorem for positive linear systems therefore gives, entrywise,
\begin{equation}\label{sd:eq:compound-comparison}
e^{-kM(t-s)}\bigwedge^kS_0(t-s)
\leq
\bigwedge^kS_Q(t,s)
\leq
e^{kM(t-s)}\bigwedge^kS_0(t-s)
\end{equation}
for $0\leq s\leq t\leq\delta$.

\textbf{Step 2. Represent the mixed determinants.}
If two endpoint samples coincide, or if two terminal indices coincide,
both determinants vanish because two rows agree. Otherwise, row
permutations affect only the sign, so order the endpoint samples as
\begin{equation*}
(n,t_1),\ldots,(n,t_s),
\qquad 0<t_1<\cdots<t_s<\delta,
\end{equation*}
and the remaining $r=n-s$ terminal samples as
\begin{equation*}
(i_1,\delta),\ldots,(i_r,\delta),
\qquad i_1<\cdots<i_r.
\end{equation*}
Set $I=(i_1,\ldots,i_r)$.  For $1\leq k\leq n$, define
\begin{equation*}
\mathcal R_k
:=(-1)^{k-1}\iota_{e_n^*}:
\bigwedge^kV\longrightarrow\bigwedge^{k-1}V.
\end{equation*}
Explicitly,
\begin{equation*}
\iota_{e_n^*}(v_1\wedge\cdots\wedge v_k)
=\sum_{j=1}^k(-1)^{j-1}e_n^*(v_j)\,
v_1\wedge\cdots\wedge\widehat{v_j}\wedge\cdots\wedge v_k.
\end{equation*}
Here the hat means that the factor is omitted.
On an ordered basis vector $e_J$, the contraction is zero if $n\notin J$.
If $n\in J$, then $n$ is the last index and
$\mathcal R_ke_J=e_{J\setminus\{n\}}$.
Consequently every $\mathcal R_k$ has nonnegative matrix entries in the
ordered exterior bases.

Set $t_0=0$, with $t_s:=t_0$ when $s=0$, and
$\mathcal T_Q^{(0)}:=\operatorname{Id}_{\bigwedge\limits^nV}$.
For $1\leq m\leq s$, define recursively
\begin{equation*}
\mathcal T_Q^{(m)}
:=\mathcal R_{n-m+1}\,
\bigwedge^{n-m+1}S_Q(t_m,t_{m-1})\,
\mathcal T_Q^{(m-1)}.
\end{equation*}
Finally, define $\mathcal T_Q:\bigwedge\limits^nV\to\bigwedge\limits^rV$ by
\begin{equation*}
\mathcal T_Q
:=\bigwedge^rS_Q(\delta,t_s)\,
\mathcal T_Q^{(s)}.
\end{equation*}
For $s=0$, this gives
$\mathcal T_Q=\bigwedge\limits^nS_Q(\delta,0)$.  With the conventions
$e_\varnothing=1$ and $\bigwedge\limits^0V=\mathbb{R}$, the definition also covers
$r=0$.

We identify the determinant by induction on the number of endpoint samples.
Before evaluating the $m$th endpoint sample, the remaining $k=n-m+1$
coordinates are propagated from $t_{m-1}$ to $t_m$ by $\bigwedge\limits^kS_Q$.
Laplace expansion along the next endpoint row applies $\iota_{e_n^*}$ to
this wedge; replacing it by $\mathcal R_k$ contributes $(-1)^{k-1}$.
After the $s$ endpoint rows, propagation to $\delta$ and extraction of the
$e_I$ component evaluate the remaining terminal rows in their stated order.
Multiplying these signs gives, with $\kappa=\sum\limits_{k=r+1}^n(k-1)$,
\begin{equation}\label{sd:eq:mixed-determinant-wedge}
\langle e_I,\mathcal T_QE_n\rangle_{\bigwedge^rV}
=(-1)^\kappa
\det\bigl[K_Q(\sigma_a,j)\bigr]_{a,j=1}^n,
\end{equation}
for the above row order.  A different order only multiplies the right-hand
side by the sign of the row permutation.

\textbf{Step 3. Compare the determinants.}
All evolution factors and contraction maps have nonnegative entries
in the ordered exterior bases. We can therefore apply
\cref{sd:eq:compound-comparison} successively to the evolution factors.
The accumulated exponent is
\begin{equation*}
\Xi
:=M\Bigl[
\sum_{m=1}^s(n-m+1)(t_m-t_{m-1})
+r(\delta-t_s)
\Bigr]
\leq nM\delta.
\end{equation*}
It follows that
\begin{equation*}
e^{-\Xi}D_0(\boldsymbol\sigma)
\leq D_Q(\boldsymbol\sigma)
\leq e^\Xi D_0(\boldsymbol\sigma).
\end{equation*}
Here \cref{sd:eq:mixed-determinant-wedge} identifies each absolute
determinant with its nonnegative matrix coefficient.
Using $\Xi\leq nM\delta$ yields \cref{sd:eq:determinant-comparison}.

The comparison argument applies to bounded measurable potentials
through the integral formulation of the evolution equations.
\end{proof}
Comparing the coefficients of the Gram determinant now gives a
quantitative observability bound.
\begin{corollary}
\label{sd:prop:robust-endpoint}
For fixed $h,n,\delta$, let $\mathfrak C_Q(\delta)$ denote the optimal constant such that
\begin{equation*}
	\lVert S_Q(\delta,0)z^0\rVert_{h,n}^2
	\leq\mathfrak C_Q(\delta)\int_0^\delta
	h|e_n^*S_Q(t,0)z^0|^2\,d t
\end{equation*}
holds for all $z^0\in\mathbb R^n$, and let $\mathfrak C_0(\delta)$ be its zero-potential counterpart. Then
\begin{equation}\label{sd:eq:robust-constant}
	\mathfrak C_Q(\delta)
	\leq ne^{4nM\delta}\mathfrak C_0(\delta).
\end{equation}
Consequently, under the length restrictions in \cref{sd:lem:unperturbed-endpoint}, there exists a constant $\mathcal C>0$, depending only on $L_-$ and $L_+$, such that every solution of $z'+A_{h,n}z=Q_n(t)z$ satisfies the observability estimate
\begin{equation}\label{sd:eq:perturbed-endpoint}
	\bigl\lVert z(\delta)\bigr\rVert_{h,n}^2
	\leq \mathcal{C}h^{-5}\exp\biggl(
	\frac{\mathcal{C}}{\delta}+4nM\delta
	\biggr)
	\int_0^\delta h\bigl\lvert z_n(t)\bigr\rvert^2\,d t.
\end{equation}
\end{corollary}
\begin{proof}
Use $W_Q$ and $H_Q$ from
\cref{sd:eq:tail-gramians,sd:eq:tail-terminal-energy}. For $\varepsilon>0$,
let $\mu_\varepsilon$ be the measure on $\Gamma_\delta$ that assigns
Lebesgue measure to the endpoint component and mass $\varepsilon$ to
each terminal point. Expanding the determinants and applying Fubini's theorem gives
Andr\'eief's identity for this mixed measure:
\begin{equation}\label{sd:eq:boundary-gram-determinant}
\det(W_Q+\varepsilon H_Q)
=\frac1{n!}\int_{\Gamma_\delta^n}
D_Q(\boldsymbol\sigma)^2
\prod_{a=1}^n\,d \mu_\varepsilon(\sigma_a).
\end{equation}
Write this determinant as $\sum\limits_{r=0}^na_{r,Q}\varepsilon^r$.
In \cref{sd:eq:boundary-gram-determinant}, the coefficient of
$\varepsilon^r$ is obtained by restricting to configurations with
exactly $r$ terminal points. Applying \cref{sd:eq:determinant-comparison}
on each such configuration gives
\begin{equation}\label{sd:eq:determinant-coefficient-comparison}
e^{-2nM\delta}a_{r,0}
\leq a_{r,Q}
\leq e^{2nM\delta}a_{r,0}.
\end{equation}
If $(z^0)^*W_Qz^0=0$, then since the observation
$t\mapsto e_n^*S_Q(t,0)z^0$ is continuous, the integral condition implies $e_n^*S_Q(t,0)z^0=0$ throughout $(0,\delta)$. 

Since $z_n(t)\equiv 0$ almost everywhere, its time derivative also vanishes. Substituting this into the $n$-th component equation of $z'+A_{h,n}z=Q_n(t)z$ forces $z_{n-1}(t)=0$ almost everywhere. By the continuity of the trajectory $z(t)$, this holds for every $t\in(0,\delta)$. Repeating this backward induction along the chain of nodes yields $z^0=0$. For $n=1$ this follows directly from the observation.
Thus $W_Q$ is positive definite and
$a_{0,Q}=\det W_Q>0$. The linear coefficient of
$\det(W_Q+\varepsilon H_Q)$ is
$a_{0,Q}\operatorname{tr}(W_Q^{-1}H_Q)$, so
\cref{sd:eq:determinant-coefficient-comparison} gives
\begin{equation*}
\operatorname{tr}(W_Q^{-1}H_Q)
=\frac{a_{1,Q}}{a_{0,Q}} 
\leq e^{4nM\delta}\operatorname{tr}(W_0^{-1}H_0) 
\leq ne^{4nM\delta}\mathfrak C_0(\delta).
\end{equation*}
By \cref{sd:eq:tail-gramians,sd:eq:tail-terminal-energy},
\begin{equation*}
\mathfrak C_Q(\delta)
=\lambda_{\max}\bigl(W_Q^{-1/2}H_QW_Q^{-1/2}\bigr)
\leq\operatorname{tr}(W_Q^{-1}H_Q),
\end{equation*}
which proves \cref{sd:eq:robust-constant}. Combining
\cref{sd:eq:robust-constant,sd:eq:unperturbed-endpoint} with $n\leq L_+/h$
gives \cref{sd:eq:perturbed-endpoint}.
\end{proof}

\bibliographystyle{plain-math}
\bibliography{references}

\end{document}